\documentclass[11pt]{article}
\usepackage{pxfonts}
\usepackage{yfonts}
\usepackage{dsfont}

\usepackage{graphicx}
\usepackage{relsize}

\def\bel{\begin{equation}\label}
\def\eeq{\end{equation}}
\def\ds{\displaystyle}
\def\endproof{\hphantom{MM}
\hfill\llap{$\square$}\goodbreak}

\def\mt{\longrightarrow}
\def\v{\vskip 1em}

\def\ve{\varepsilon}
\def\R{\mathds R}
\def\Z{\mathds Z}
\def\C{\mathfrak{B}}
\def\Cc{\mathfrak{C}}

\def\N{{\bf N}}

\def\exp{{\bf exp}}

\def\S{{\bf S}}

\def\Sz{\mathfrak{S}}

\def\E{{\bf E}}

\def\P{{\bf P}}
\def\Q{{\bf Q}}
\def\D{\mathfrak{D}}

\def\J{{\bf J}}

\def\B{{\bf B}}
\def\H{{\bf H}}
\def\L{{\bf L}}
\def\U{{\bf U}}
\def\V{{\bf V}}
\def\T{{\bf T}}

\def\p{{\partial}}

\def\i{{\bf i}}
\def\Tilde{\widetilde}
\def\Hat{\widehat}

\def\bar{\overline}
\def\supp{{\bf supp}}

\def\I{{\bf I}}

\def\M{{\bf M}}

\def\alpha{\alphaup}
\def\beta{\betaup}
\def\gamma{\gammaup}
\def\delta{\deltaup}
\def\xi{{\xiup}}
\def\eta{{\etaup}}
\def\tau{{\tauup}}
\def\rho{{\rhoup}}
\def\phi{{\phiup}}
\def\psi{{\psiup}}
\def\lambda{{\lambdaup}}
\def\omega{\omegaup}
\def\varphi{{\varphiup}}
\def\gamma{{\gammaup}}
\def\c{{\bf c}}
\def\t{{\bf t}}
\def\s{{\bf s}}

\def\a{{\bf a}}
\def\b{{\bf b}}

\def\p{{\partial}}
\def\vv{{\bf v}}
\def\h{{\bf h}}
\def\({\left(}
\def\){\right)}

\newtheorem{cor}{Corollary}[section]
\newtheorem{lemma}{Lemma}[section]

\newtheorem{remark}{Remark}[section]

\begin{document}
\[\begin{array}{cc}\ds\hbox{\LARGE{\bf Singular integrals of non-convolution type} }
  \\\\ \ds
 \hbox{\LARGE{\bf on product spaces}}
 \end{array}\]

\[\hbox{Zipeng Wang}\]
\begin{abstract}
We study a new class of pseudo differential operators whose symbols satisfy the differential inequality with a mixture of homogeneities. On the other hand, by taking singular integral realization, it can be equivalently defined by kernels carrying certain characteristic properties on product spaces. We prove that these operators are bounded on $\L^p$-spaces  for $1<p<\infty$. Moreover, they form an algebra.
\end{abstract}

\section{Introduction}
\setcounter{equation}{0}
The study of certain operators  commuting with a multi-parameter family of dilations  dates back to the time of  Jessen, Marcinkiewicz and Zygmund.  Over the several past  decades,
 a number of pioneering  results  have been accomplished, for example   
by  Fefferman \cite{R.Fefferman}-\cite{R.Fefferman''}, Fefferman and Stein \cite{R-F.S},  Cordoba and Fefferman \cite{Cordoba-Fefferman}, Chang and  Fefferman \cite{Chang-Fefferman}, Journ\'{e} \cite{Journe'},  Pipher \cite{Pipher}, Carbery and Seeger \cite{C.S}-\cite{C.S'}, Nagel and Stein \cite{N.S} and M\"{u}ller, Ricci and Stein \cite{M.R.S}. 

In this paper, we introduce a new class of pseudo differential operators  defined either by symbols or kernels alternatively. 

Let $\R^\N=\R^{\N_1}\times\R^{\N_2}\times\cdots\times\R^{\N_n}$ for $n\ge2$. Consider
\bel{T}
\begin{array}{lr}\ds
\T f(x)~=~\int_{\R^\N} e^{2\pi\i x\cdot\xi} \sigma(x,\xi)\Hat{f}(\xi)d\xi
\end{array}
\eeq
where $\sigma(x,\xi)\in\mathcal{C}^\infty(\R^\N\times\R^\N)$. 

$\diamond$ {\small Throughout,  $\C>0$ is a  generic constant whose value depends on the subindices.}

{\bf Symbol class $\S_\rhoup$:}  ~~We say  $\sigma\in\S_\rhoup, 0<\rho<1$ if 
\bel{Ineq}
\left|\p_\xi^\alphaup\p_x^\betaup\sigma(x,\xi)\right|~\leq~\C_{\alphaup~\betaup}~\prod_{i=1}^n \Bigg({1\over 1+|\xi_i|+|\xi|^{ \rhoup}}\Bigg)^{|\alphaup_i|}(1+|\xi|)^{\rhoup|\betaup|}
\eeq
for every multi-indices $\alphaup, \betaup$. 

The original version of $\S_\rho$ for $\rho={1\over 2}$ is first suggested by Stein. This symbol class arises by considering compositions of Calder\'{o}n-Zygmund operators with different homogeneities of given dilations, i.e: one isotropic and the other parabolic. A systematic discussion for  certain classes of operators  having  multi-parameter characteristics has been established in the Memoirs paper by Nagel, Ricci, Stein and Wainger \cite{N.R.S.W}.

Observe that $\S_\rho\subset\S_{\rho,\rho}$ so-called  the {\it exotic} class: a symbol $\sigma\in\S_{\rho,\rho}$ satisfies
\bel{Ineq exotic}
\left|\p_\xi^\alphaup\p_x^\betaup\sigma(x,\xi)\right|~\leq~\C_{\alphaup~\betaup}~(1+|\xi|)^{-\rhoup|\alpha|+\rhoup|\betaup|}
\eeq
for every multi-indices $\alphaup, \betaup$. 
In particular, the Fourier transform of  $\hbox{\bf R}$: $\hbox{\bf R}(x)=e^{\i/|x|}\left({1\over|x|}\right)^{{3\over 2}}, x\neq0$  belongs $\S_{{1\over 2},{1\over 2}}$.  This distribution is known as {\it Riemann singularity}. Convolution with $\hbox{\bf R}$ is  bounded  only on $\L^2$-spaces. Regarding details can be found in chapter  VII of  Stein \cite{S}.

From (\ref{T}), we rewrite 
\bel{T Omega}
\begin{array}{cc}\ds
\T f(x)~=~\int_{\R^\N} f(y) \Omega(x,y)dy,
\\\\ \ds
\Omega(x,y)~=~\int_{\R^\N} e^{2\pi\i (x-y)\cdot\xi}\sigma(x,\xi)d\xi,\qquad x\neq y.
\end{array}
\eeq
 \begin{remark} Denote $z=x-y$ and  $\Omega(x,y)=\Omega(x,z)$. We have $\sigma\in\S_\rho, 0<\rho<1$ satisfying (\ref{Ineq})  equivalent to $\Omega$ having the size estimate and cancellation property given in below.
\end{remark}

{\bf Size estimate} ~~{\it For every multi-index $\alphaup$ and $\betaup$,  we have
 \bel{Omega Est}
\left|\p_z^\alphaup\p_x^\betaup\Omega(x,z)\right|~\leq~\C_{\alpha~\beta}~\left({1\over |z|}\right)^{\rho|\beta|}\prod_{i=1}^n \Bigg({1\over |z_i|+|z|^{1/\rho}}\Bigg)^{\N_i+|\alpha_i|},\qquad |z|>0.
\eeq
On the other hand, $\p_z^\alphaup\p_x^\betaup\Omega(x,z)$ decays rapidly  as $|z|\mt\infty$. }

Let  $\varphi_i\in\mathcal{C}^\infty_o(\R^{\N_i}), i=1,2,\ldots,n$   be normalized {\it bump}-functions.

{\bf Cancellation property} ~~{\it We have
\bel{Cancellation}
\begin{array}{lr}\ds
\left|\int _{\R^\N}\p_z^\alpha\p_x^\betaup\Omega(x,z) \prod_{i=1}^n \varphi(R_iz_i) dz\right|
~\leq~\C_{\alphaup~\betaup}~\left\{ 1+\sum_{i=1}^n R_i\right\}^{\rhoup|\betaup|+|\alpha|}
\end{array}
\eeq
for every multi-indices $\alphaup$, $\betaup$ and  every $R_i>0, i=1,2,\ldots,n$.}

{\bf Theorem One} ~ {\it Let $\T$ either defined in (\ref{T}) or (\ref{T Omega}). Suppose  $\sigma\in\S_\rho, 0<\rho<1$ satisfying (\ref{Ineq}) or equivalently $\Omega$ satisfies (\ref{Omega Est})-(\ref{Cancellation}).  We have
\bel{Result One}
\left\| \T f\right\|_{\L^p(\R^\N)}~\leq~\C_{p~\rhoup}~\left\| f\right\|_{\L^p(\R^\N)},\qquad 1<p<\infty.
\eeq
Moreover, this class of operators form an algebra under composition.}

Our paper is organized as follows.  In the next section, we introduce a framework of Littlewood-Paley projections and  obtain certain combinatorial estimates on the regarding multi-parameter dyadic decomposition. In section 3, we give a classification between  $\sigma\in\S_\rhoup$ and  $\Omega$ satisfying (\ref{Omega Est})-(\ref{Cancellation}).  
In section 4, we prove  a fundamental lemma.  As a corollary,  $\T$ with $\sigma\in\S_\rho$ form an algebra.
In section 5, we show that every  partial operator is bounded by the strong maximal function. Furthermore, these operators enjoy a desired property of almost orthogonality. 
In section 6, we conclude the $\L^p$-norm inequality in (\ref{Result One}). 
In section 7,  we prove the required Littlewood-Paley inequality.

\section{Combinatorial estimates on the regarding dyadic decomposition}
\setcounter{equation}{0}
Let  $t_i, i=1,2,\ldots,n$ be integers. Denote $\t=(t_1,t_2,\ldots,t_n)$. For $0<\rho<1$, we write $q=1/\rho$  and
\bel{t_i dila}
\begin{array}{rl}\ds
\t_i\xi~=~\left(2^{-q t_i}\xi_1,~\ldots,~2^{-t_i}\xi_i,~\ldots,~2^{-q t_i}\xi_n\right),
~~~~
i=1,2,\ldots,n.
\end{array}
\eeq
Consider  $\varphi\in\mathcal{C}^\infty_o(\R)$ such that $\varphi(s)=1$ if $|s|\leq1$ and $\varphi(s)=0$ for $|s|>2$. 
 We define
\bel{delta_t}
\deltaup_\t(\xi)~=~ \prod_{i=1}^n \phi\left(\t_i\xi\right),\qquad \phi(\xi)~=~\varphi\left(|\xi|\right)-\varphi\left(2|\xi|\right).
\eeq 
The support of $\deltaup_\t(\xi)$ is contained in the intersection of $n$ elliptical shells, with different homogeneities of given dilations. 

A partial  operator $\Delta_\t$ is defined by
\bel{Delta_t}
\Hat{\Delta_{\t} f}(\xi)~=~ \deltaup_\t\left(\xi\right)\Hat{f}(\xi).
\eeq
Let 
\bel{t_imath}
t_\imath~=~\max\Big\{t_i\colon i=1,2,\ldots,n\Big\}.
\eeq
We partition the set $\{1,2,\ldots,n\}=\I\cup\J$  for which
\bel{IJ}
t_\imath~\leq~ q t_i-\left(2+\log_2n\right),\qquad i\in\I\qquad\hbox{and}\qquad  t_\imath~>~q t_i-\left(2+\log_2n\right),\qquad i\in\J.
\eeq
{\bf Index class H:} ~~{\it We say $\t\in\H$ if there exists at least one $i\in\{1,2,\ldots n\}$ such that
\bel{H} 
t_i~\ge~ {1\over q-1} \left(2+\log_2n\right).
\eeq}
\begin{lemma}  
  Let $\xi\in\supp\deltaup_\t(\xi)$. Suppose $\t\in\H$. We have
\bel{Combinatorial results}
\begin{array}{cc}\ds
|\xi_i|~\sim~2^{t_i},\qquad i\in\I,
\\\\ \ds
|\xi_i|~\lesssim~2^{t_i},\qquad|\xi_\imath|~\sim~2^{t_\imath}~\sim~2^{q t_i},\qquad i\in\J.
\end{array}
\eeq
\end{lemma}

{\bf Proof :}  Note that  $\t\in\H$ satisfying (\ref{H}) guarantees that the set $\I$ is non-empty. In particular, $\imath\in\I$ because of (\ref{t_imath}).

Consider $\xi\in\supp\deltaup_\t(\xi)$. We write
$\xi=\left(\xi_i,\xi_i^\dagger\right)\in\R^{\N_i}\times\R^{\N-\N_i}$ for $i=1,2,\ldots,n$.

 By definition of $\deltaup_\t(\xi)$ in (\ref{delta_t}), we find  
\bel{xi_i est in supp}
\begin{array}{lr}\ds
 |\xi_i|~<~2^{t_i+1},\qquad i=1,2,\ldots,n\qquad\hbox{and}
 \\\\ \ds
 |\xi_i|~<~2^{q t_j+1},\qquad\hbox{for}~~ j\neq i.
 \end{array}
 \eeq 
On the other hand,  we either have
\bel{est Dila 1}
|\xi_i|~>~{2^{t_i-1}\over \sqrt{2}}\qquad\hbox{or}\qquad \left|\xi_i^\dagger\right|~>~{2^{q t_i-1}\over \sqrt{2}},\qquad i=1,2,\ldots,n.
\eeq

Let $i\in\I$ in (\ref{IJ}). Suppose  $|\xi_i|\leq {2^{t_i-1}\over \sqrt{2}}$. From (\ref{est Dila 1}),
there is some $\xi_j$ for  $j\neq i$ such that 
\bel{est Dila 1.1}
|\xi_j|~>~{2^{q t_i-1}\over \sqrt{2}}{1\over \sqrt{n-1}}.
\eeq
Together with $|\xi_j|<2^{t_j+1}$ shown in (\ref{xi_i est in supp}), we must have
\bel{ineq case 1}
q t_i-2-{1\over 2}\log_2 2(n-1)~<~t_j.
\eeq
However,  $i\in\I$ implies 
$t_j\leq t_\imath \leq q t_i-\left(2+\log_2n\right)$.
The inequality in (\ref{ineq case 1}) cannot be true because 
\bel{fact recall}
(1/2)\log_2 2(n-1)~<~\log_2 n\qquad\hbox{for}\qquad n\ge2.
\eeq
Therefore, by putting together all estimates above, we find
\bel{xi_i I norm}
{2^{t_i-1}\over \sqrt{2}}~<~|\xi_i|~<~2^{t_i+1},\qquad i\in\I.
\eeq
Recall $\imath\in\I$.  From (\ref{xi_i I norm}) and the second inequality in (\ref{xi_i est in supp}), we find $ 2^{t_\imath}\sim|\xi_\imath|<2^{q t_i+1}$ for $i\neq \imath$.   Consider $i\in\J$ in (\ref{IJ}) where $q t_i-\left(2+\log_2n\right)<t_\imath$. We necessarily have
\bel{i in J range}
2^{t_\imath}~\sim~|\xi_\imath|~\sim~2^{q t_i},\qquad i\in\J.
\eeq
\endproof
\begin{lemma}
Let $\sigma\in\S_\rhoup,0<\rho<1$. Suppose $\t\in\H$. We have 
\bel{Compo est symbol}
\begin{array}{lr}\ds
\left| \p_\xi^\alphaup\deltaup_\t(\xi)\sigma(x,\xi)\right|~\leq~\C_\alphaup ~\prod_{i=1}^n\left({1\over 1+|\xi_i|+|\xi|^\rhoup}\right)^{|\alphaup_i|}
~\sim~\C_\alphaup ~\prod_{i=1}^n 2^{-t_i|\alphaup_i|}
\end{array}
\eeq
for every multi-index $\alphaup$.
\end{lemma}
{\bf Proof:} First, $\deltaup_\t(\xi)$ defined in (\ref{delta_t}) clearly satisfies 
\bel{Compo est symbol delta}
\begin{array}{lr}\ds
\left| \p_\xi^\alphaup\deltaup_\t(\xi)\right|~\leq~\C_\alphaup ~\prod_{i=1}^n 2^{-t_i|\alphaup_i|}
\end{array}
\eeq
for every multi-index $\alphaup$.

Let $q=1/\rho$,  $t_\imath=\max_{i\in\{1,2,\ldots,n\}} t_i$ and $\I\cup\J=\{1,2,\ldots,n\}$ defined in (\ref{IJ}).  

Recall {\bf Lemma 2.1}. We  have
\bel{Lemma 2.1 recall}
\begin{array}{lr}\ds
|\xi_i|~\sim~ 2^{t_i},\qquad |\xi|~\sim~|\xi_\imath|~\sim~2^{t_\imath}~\lesssim~2^{q t_i},\qquad i\in\I,
\\\\ \ds
|\xi_i|~\lesssim~2^{t_i},\qquad |\xi|~\sim~|\xi_\imath|~\sim~2^{t_\imath}~\sim~ 2^{qt_i},\qquad  i\in\J
\end{array}
\eeq
provided that $\xi\in \supp \deltaup_\t(\xi)\sigma(x,\xi)$ and $\t\in\H$.

From (\ref{Lemma 2.1 recall}), we conclude
\bel{condition i}
1+|\xi_i|+|\xi|^\rho~\sim~ 2^{t_i},\qquad i=1,2,\ldots,n.
\eeq 
\endproof

\section{Classification for symbols and kernels}
\setcounter{equation}{0}
Let $\sigma\in\S_\rhoup, 0<\rho<1$ satisfying the differential inequality in (\ref{Ineq}). Consider
\bel{Omega z}
\Omega(x,z)~=~\int_{\R^\N} e^{2\pi\i z\cdot\xi} \sigma(x,\xi) d\xi.
\eeq
For convention, we fix some notations as follows.

$\bullet$ Given $\xi_i\in\R^{\N_i}$ and a multi-index $\alpha_i$, we denote $\xi_i^{\alpha_i}=\prod_{j=1}^{\N_i} \xi_{i j}^{\alpha_{ij}}$ for every $i=1,2,\ldots,n$. Moreover, $\xi^\alpha=\prod_{i=1}^n \xi_i^{\alpha_i}$.

$\bullet$ Let $\imath\in\{1,2,\ldots,n\}$ such that $|z_\imath|=\max_{i=1,2,\ldots,n}|z_i|$. 

$\bullet$ For each $i=1,2,\ldots,n$, $|z_{i\jmath}|= \max_{j=1,2,\ldots,\N_i}|z_{ij}|$ and $|\xi_{i\jmath}|= \max_{j=1,2,\ldots,\N_i}|\xi_{ij}|$.

$\bullet$ $\gamma=(\gamma_1,\gamma_2,\ldots,\gamma_n)$ is a multi-index where $\gamma_i=\gamma_{i\jmath}\ge0, i=1,2,\ldots,n$. 

$\bullet$ In particular,  we write $\gamma_\imath=\gamma_{\imath\jmath}=\lambda_1+\lambda_2+\cdots+\lambda_n$ for which  $z_\imath^{\lambda_i}=z_{\imath\jmath}^{\lambda_i}$ and $\p_{\xi_\imath}^{\lambda_i}=\p_{\xi_{\imath\jmath}}^{\lambda_i}$ where $\lambda_i\ge0, i=1,2,\ldots,n$.

Recall $\deltaup_\t(x)$ defined in (\ref{delta_t}) and  the index class {\bf H} from (\ref{H}). 
Consider
\bel{partial Omega}
 \begin{array}{lr}\ds
\p_z^\alphaup\p_x^\betaup\Omega(x,z)
~=~ \left(2\pi\i\right)^{|\alphaup|}\int_{\R^\N} e^{2\pi\i z\cdot\xi}  \xi^\alphaup  \p_x^\betaup\sigma(x,\xi)d\xi
\\\\ \ds~~~~~~~~~~~~~~~~~~~
~=~ \left( 2\pi\i\right)^{|\alphaup|}\sum_{\t\in\H}   \p_z^\alphaup\p_x^\betaup\Omega_\t(x,z)~+ ~
\left( 2\pi\i\right)^{|\alphaup|}\p_z^\alphaup\p_x^\betaup\Omega^\flat(x,z)
 \end{array}
 \eeq
 where 
\bel{partial Omega_t}
\p_z^\alphaup\p_x^\betaup\Omega_\t(x,z)~=~\int_{\R^\N} e^{2\pi\i z\cdot\xi} \xi^\alphaup  \deltaup_\t(\xi)\p_x^\betaup\sigma(x,\xi)d\xi
\eeq
and
\bel{Omega flat}
 \p_z^\alphaup\p_x^\betaup\Omega^\flat(x,z)~=~ \left({1\over 2\pi\i}\right)^{-|\alphaup|}   \int_{\R^\N} e^{2\pi\i z\cdot\xi} \xi^\alphaup \Tilde{\deltaup}(\xi)\p_x^\betaup\sigma(x,\xi)d\xi,~~~~  \Tilde{\deltaup}(\xi)=\sum_{\t\notin\H} \delta_\t(\xi).
 \eeq 
First, we focus on the summation  of $\t\in\H$ in (\ref{partial Omega}).

Let $\gamma$ be the multi-index described as above. Consider
\bel{Kernel size rewrite}
\begin{array}{lr}\ds
\p_z^\alphaup\p_x^\betaup\Omega_\t(x,z)
~=~\prod_{i=1}^n z_i^{-\gamma_i} \V^{\alpha~\beta~\gamma}_\t(z)
\end{array}
\eeq
where
\bel{V_t}
\begin{array}{lr}\ds
\V^{\alpha~\beta~\gamma}_\t(z)~=~
 \prod_{i=1}^n z_i^{\gamma_i}\int_{\R^\N} e^{2\pi\i z\cdot\xi}  \xi^\alphaup  \deltaup_\t(\xi)\p_x^\betaup\sigma(x,\xi)d\xi
\\\\ \ds~~~~~~~~~~~~~~
~=~z_\imath^{\lambda_\imath}\prod_{i\neq\imath} z_i^{\gamma_i}z_\imath^{\lambda_i}\int_{\R^\N} e^{2\pi\i z\cdot\xi} \xi_\imath^{\alpha_\imath}\prod_{i\neq \imath} \xi_i^{\alphaup_i}  \deltaup_\t(\xi)\p_x^\betaup\sigma(x,\xi)d\xi.
\end{array}
\eeq
By integration by parts $w.r.t~ \xi$, we find 
\bel{Kernel size Int by parts}
\begin{array}{lr}\ds
\V^{\alpha~\beta~\gamma}_\t(z)~=~\left({-1\over 2\pi\i}\right)^{|\gamma|}
\int_{\R^\N} e^{2\pi\i z\cdot\xi} \p_{\xi_\imath}^{\lambda_\imath}\xi_\imath^{\alpha_\imath}\prod_{i\neq \imath}\p_{\xi_i}^{\gamma_i}\p_{\xi_\imath}^{\lambda_i} \xi_i^{\alphaup_i}  \deltaup_\t(\xi)\p_x^\betaup\sigma(x,\xi)d\xi.
\end{array}
\eeq
Recall that $\sigma(x,\xi)$ satisfies the differential inequality in (\ref{Ineq}).
By using {\bf Lemma 2.2}.  we find  
\bel{combin est'} 
\begin{array}{lr}\ds
\left| \V^{\alpha~\beta~\gamma}_\t(z)\right|~\leq~\C_{\alpha~\beta~\gamma}~
2^{-\lambda_\imath t_\imath+\left(\N_\imath+|\alpha_\imath|+\rho|\beta|\right)t_\imath}\prod_{i\neq\imath} 2^{-\gammaup_i t_i-\lambda_i t_\imath+\left(\N_i +|\alphaup_i|\right)t_i}.
\end{array}
\eeq
From {\bf Lemma 2.1}, we have $2^{t_\imath}\lesssim 2^{t_i/\rhoup}$ for every $i=1,2,\ldots,n$. 
Therefore, (\ref{combin est'}) can be further bounded by
\bel{combin est}
\C_{\alpha~\beta~\gamma}~
2^{-\lambda_\imath t_\imath+\left(\N_\imath+|\alpha_\imath|+\rho|\beta|\right)t_\imath}\prod_{i\neq\imath} 2^{-\left(\gammaup_i+\rho\lambda_i\right)t_i+\left(\N_i +|\alphaup_i|\right)t_i}.
\eeq
Choose $\lambda_\imath=0$ if
 $2^{t_\imath}\leq |z_\imath|^{-1}$ and $\lambda_\imath>\N_\imath+|\alphaup_\imath|+\rhoup|\betaup|$ if $2^{t_\imath}> |z_\imath|^{-1}$. We have
 \bel{Sum imath}
 \begin{array}{lr}\ds
|z_\imath|^{-\lambda_\imath} \sum_{t_\imath\ge0} 2^{-\lambda_\imath t_\imath+\left(\N_\imath+|\alpha_\imath|+\rho|\beta|
\right)t_\imath}
\\\\ \ds
~\leq~\sum_{2^{t_\imath}\leq |z_\imath|^{-1}} 2^{\left(\N_\imath+|\alpha_\imath|+\rho|\beta|
\right)t_\imath}~+~|z_\imath|^{-\lambda_\imath}\sum_{2^{t_\imath}>|z_\imath|^{-1}} 2^{-\lambda_\imath t_\imath+\left(\N_\imath+|\alpha_\imath|+\rho|\beta|
\right)t_\imath}
\\\\ \ds
~\leq~\C~\left({1\over |z_\imath|}\right)^{\N_\imath+|\alpha_\imath|+\rho|\beta|}.
 \end{array}
 \eeq 
Let $\lambda_i=0$ for $i\neq\imath$. Choose $\gammaup_i=0$ for $2^{t_i}\leq |z_i|^{-1}$ and $\gammaup_i>\N_i+|\alphaup_i|$ for $2^{t_i}> |z_i|^{-1}$. 
We have
 \bel{Sum i}
 \begin{array}{lr}\ds
|z_i|^{-\gammaup_i} \sum_{t_i\ge0} 2^{-\gammaup_i t_i+\left(\N_i+|\alpha_i|\right)t_i}
\\\\ \ds
~\leq~\sum_{2^{t_i}\leq |z_i|^{-1}} 2^{\left(\N_i+|\alpha_i|
\right)t_i}~+~|z_i|^{-\gammaup_i}\sum_{2^{t_i}>|z_i|^{-1}} 2^{-\gammaup_i t_i+\left(\N_i+|\alpha_i|
\right)t_i}
\\\\ \ds
~\leq~\C~\left({1\over |z_i|}\right)^{\N_i+|\alpha_i|}.
 \end{array}
 \eeq 
Let $\gammaup_i=0$ for every $i\neq\imath$. Choose $\lambda_i=0$ for $2^{t_i}\leq |z_\imath|^{-{1/\rhoup}}$ and 
$\lambda_i>\left(\N_i+|\alphaup_i|\right)/\rhoup$ for $2^{t_i}> |z_\imath|^{-{1/\rhoup}}$. We have
 \bel{Sum i rho}
 \begin{array}{lr}\ds
|z_\imath|^{-\lambda_i} \sum_{t_i\ge0} 2^{-\rho\lambda_i t_i+\left(\N_i+|\alpha_i|\right)t_i}
\\\\ \ds
~\leq~\sum_{2^{t_i}\leq |z_\imath|^{-{1/ \rho}}} 2^{\left(\N_i+|\alpha_i|
\right)t_i}~+~|z_\imath|^{-\lambda_i}\sum_{2^{t_i}>|z_\imath|^{-{1/ \rho}}} 2^{-\rho \lambda_i t_i+\left(\N_i+|\alpha_i|
\right)t_i}
\\\\ \ds
~\leq~\C~\Bigg({1\over |z_\imath|^{1/\rho}}\Bigg)^{\N_i+|\alpha_i|}.
 \end{array}
 \eeq 
Recall (\ref{Kernel size rewrite})-(\ref{Kernel size Int by parts}). By putting together the above estimates from (\ref{combin est'}) to (\ref{Sum i rho}), we find
\bel{Size Est'}
\begin{array}{lr}\ds
 \sum_{\t\in\H} \left| \p_z^\alphaup\p_x^\betaup\Omega_\t(x,z)\right|~\leq~\C_{\alpha~\beta}~\left({1\over |z_\imath|}\right)^{\N_\imath+|\alpha_\imath|+\rho|\beta|}\prod_{i\neq\imath} \Bigg({1\over |z_i|+|z_\imath|^{1/\rho}}\Bigg)^{\N_i+|\alpha_i|}
\end{array}
\eeq
for $|z|>0$. 
 Note that $|z|\sim|z_\imath|=\max_{i\in\{1,2,\ldots,n\}}|z_i|$. 
 
 From (\ref{Size Est'}), we obtain
\bel{Size Est}
\begin{array}{lr}\ds
 \sum_{\t\in\H} \left| \p_z^\alphaup\p_x^\betaup\Omega_\t(x,z)\right|~\leq~\C_{\alpha~\beta}~\left({1\over |z|}\right)^{\rho|\beta|}\prod_{i=1}^n \Bigg({1\over |z_i|+|z|^{1/\rho}}\Bigg)^{\N_i+|\alpha_i|},\qquad |z|>0.
\end{array}
\eeq
On the other hand, by integration by parts $w.r.t~\xi$, we find that
$\sum_{\t\in\H}  \p_z^\alphaup\p_x^\betaup\Omega_\t(x,z)$ decays rapidly as $|z|\mt\infty$. 

Return to $\Omega^\flat(x,z)$ defined in (\ref{Omega flat}). Note that $\Tilde{\delta}(\xi)=\sum_{\t\notin\H}\deltaup_\t(\xi)$ has a compact support depending on the index class {\bf H}  given in (\ref{H}). 
By definition of $\delta_\t(\xi)$ in (\ref{delta_t}), we find   $|\xi|\lesssim 2^{\rhoup\over 1-\rhoup}$ for $0<\rhoup<1$ whenever $\xi\in\supp\Tilde{\deltaup}(\xi)$. 
This implies
\bel{Omega flat norm bound}
 \left| \p_z^\alpha \p_x^\beta \Omega^\flat(x,z)\right|~\leq~\C_{\alpha~\beta}~2^{{\rho\over 1-\rho}\left[\N+|\alpha|+\rho|\beta|\right]}.
 \eeq
From (\ref{t_i dila})-(\ref{delta_t}), we have
\bel{Tilde delta rewrite}
\begin{array}{lr}\ds
\Tilde{\delta}(\xi)~=~\sum_{\t\notin\H}\deltaup_\t(\xi)
\\\\ \ds~~~~~~
~=~\prod_{i=1}^n\varphi\left(\sqrt{2^{-2 t/\rho}|\xi_1|^2+\cdots+2^{-2t}|\xi_i|^2+\cdots+2^{-2 t/\rho}|\xi_n|^2}\right)
\end{array}
\eeq
where $\varphi\in\mathcal{C}^\infty_o(\R)$ is the smooth {\it bump}-function as before and $t$ is the largest integer such that $t<{\rho\over 1-\rho}\left(2+\log_2 n\right)$.
Every derivative of $\Tilde{\delta}(\xi)$ is either bounded by $\C 2^{-t}$ when $|\xi|\sim2^t$ or zero otherwise.

Let $\sigma\in\S_\rho$ satisfying  (\ref{Ineq}). Rewrite
\bel{Omega flat rewrite}
\begin{array}{lr}\ds
 \p_z^\alphaup\p_x^\betaup\Omega^\flat(x,z)~=~ \left({1\over 2\pi\i}\right)^{-|\alphaup|} \sum_{j\in\Z}  \int_{\R^\N} e^{2\pi\i z\cdot\xi} \xi^\alphaup \Tilde{\deltaup}(\xi)\p_x^\betaup\sigma(x,\xi)\varphi\left(2^{-j}|\xi|\right)d\xi.
 \end{array}
 \eeq 
Recall $|z_\imath|=\max_{i\in\{1,2,\ldots,n\}}|z_i|$. An $N$-fold integration by parts $w.r.t~\xi_\imath$ inside each integral above shows
\bel{flat int est}
\begin{array}{lr}\ds
\left|\int_{\R^\N} e^{2\pi\i z\cdot\xi} \xi^\alphaup \Tilde{\deltaup}(\xi)\p_x^\betaup\sigma(x,\xi)\varphi\left(2^{-j}|\xi|\right)d\xi\right|
~\leq~\C_N~\left({1\over |z|}\right)^N~2^{j\left[\N+|\alpha|+\rho|\beta|\right]}2^{-j\rho N}.
\end{array}
\eeq
Choose $N=0$ for $2^j\leq|z|^{-1/\rho}$ and $N>(\N+|\alpha|+\rho|\beta|)/\rho$ for $2^j>|z|^{-1/\rho}$. We have
 \bel{Sum i rho flat}
 \begin{array}{lr}\ds
 \left|  \p_z^\alphaup\p_x^\betaup\Omega^\flat(x,z)\right|~\leq~\C_{\alpha~N}  \left({1\over |z|}\right)^N~2^{j\left[\N+|\alpha|+\rho|\beta|\right]}2^{-j\rho N} 
 \\\\ \ds
~\leq~\C_\alpha~\sum_{2^j\leq |z|^{-{1/ \rho}}} 2^{j\left[\N+|\alpha|+\rho|\beta|
\right]}~+~ \C_{\alpha~N}\left({1\over |z|}\right)^N \sum_{2^j>|z|^{-{1/ \rho}}} 2^{j\left[\N+|\alpha|+\rho|\beta|\right]}2^{-j\rho N}
\\\\ \ds
~\leq~\C_{\alpha~\beta}~\Bigg({1\over |z|^{1/\rho}}\Bigg)^{\N+|\alpha|+\rho|\beta|}.
 \end{array}
 \eeq 
Together with  (\ref{Omega flat norm bound}), $\p_z^\alphaup\p_x^\betaup\Omega^\flat(x,z)$ satisfies the norm inequality in (\ref{Sum i rho flat}) for $|z|^{-1}\lesssim 2^{\rho^2\over 1-\rho}$. Compare (\ref{Sum i rho flat}) to (\ref{Size Est}). From direct computation, we find
\bel{Compara flat}
 \left|  \p_z^\alphaup\p_x^\betaup\Omega^\flat(x,z)\right|~\leq~\C_{\alpha~\beta}~2^{\rho\left[\N+|\alpha|+\rho|\beta|\right]}\left({1\over |z|}\right)^{\rho|\beta|}\prod_{i=1}^n \Bigg({1\over |z_i|+|z|^{1/\rho}}\Bigg)^{\N_i+|\alpha_i|},~~~~|z|>0.
 \eeq
From (\ref{partial Omega}) to (\ref{Compara flat}), we conclude (\ref{Omega Est}).

Let $\varphi_i\in\mathcal{C}^\infty_o(\R^{\N_i})$, $i=1,2,\ldots,n$ be normalized  {\it bump}-functions. 
 We consider
\bel{zeta}
\begin{array}{lr}\ds
\zeta(x,R_1,\ldots,R_n)~=~\int_{\R^\N} \p_z^\alpha\p_x^\betaup\Omega(x,z)\prod_{i=1}^n\varphi_i\left(R_iz_i\right)dz_i
\\\\ \ds ~~~~~~
~=~(-2\pi\i)^{|\alpha|}\int_{\R^\N}\p_x^\beta \sigma(x, \xi)  \prod_{i=1}^n R_i^{-\N_i}\Hat{\varphi}_i\left(-R^{-1}_i\xi_i\right)\xi_i^{\alpha_i}d\xi_i.
\end{array}
\eeq
Define 
\bel{norm function vartheta}
\vartheta_{\alpha~\betaup}(\xi)~=~1+\sum_{i=1}^n |\xi_i|^{\rhoup|\betaup|+|\alpha|}.
\eeq
Rewrite $\zeta(x,R_1,\ldots,R_n)$ in (\ref{zeta}) as
\bel{est Cancellation}
\begin{array}{lr}\ds
(-2\pi\i)^{|\alpha|}\sum_{i=1}^n \int_{\R^\N}|\xi_i|^{\rhoup|\betaup|+|\alpha|} \left\{\prod_{i=1}^n\xi_i^{\alpha_i}\p_x^\betaup\sigma(x, \xi)\vartheta^{-1}_{\alpha~\betaup}(\xi)\right\} \prod_{i=1}^nR_i^{-\N_i}\Hat{\varphi}_i\left(-R^{-1}_i\xi_i\right)d\xi
\\\\ \ds
~+~(-2\pi\i)^{|\alpha|}\int_{\R^\N} \left\{\prod_{i=1}^n\xi_i^{\alpha_i}\p_x^\betaup\sigma(x, \xi)\vartheta^{-1}_{\alpha~\betaup}(\xi)\right\} \prod_{i=1}^nR_i^{-\N_i}\Hat{\varphi}_i\left(-R^{-1}_i\xi_i\right)d\xi.
\end{array}
\eeq
By changing dilations $\xi_i\mt R_i\xi_i$ inside (\ref{est Cancellation}), we find 
\bel{zeta norm}
\begin{array}{lr}\ds
\left|\zeta(x,R_1,\ldots,R_n)\right|~\leq~\C_{\alpha~\beta}~\left\{1+\sum_{i=1}^nR_i\right\}^{\rhoup|\betaup|+|\alpha|}.
\end{array}
\eeq
From (\ref{zeta}) to (\ref{zeta norm}),
we conclude (\ref{Cancellation}).

Suppose that  $\Omega(x,z)$ satisfies   (\ref{Omega Est}) and (\ref{Cancellation}). 
 Define
\bel{rho_i}
\rho_i(\xi)~=~1+|\xi_i|+|\xi|^\rhoup,\qquad  i=1,2,\ldots,n.
\eeq
Consider
\bel{Est symbol}
\begin{array}{lr}\ds
\p_\xi^\alpha\p_x^\beta\sigma(x,\xi)\prod_{i=1}^n \left(\rho_i(\xi)\right)^{|\alpha_i|}
~=~(-2\pi\i)^{|\alphaup|}\int_{\R^\N}  \p_x^\betaup \Omega(x,z) e^{-2\pi\i z\cdot\xi}\prod_{i=1}^n  \left(\rho_i(\xi)\right)^{|\alpha_i|}z_i^{\alpha_i} dz.
\end{array}
\eeq
We aim to show 
\bel{symbol result}
\left| \p_\xi^\alpha\p_x^\beta\sigma(x,\xi)\prod_{i=1}^n \left(\rho_i(\xi)\right)^{|\alpha_i|}\right|~\leq~\C_{\alpha~\beta}~\left(1+|\xi|\right)^{\rho|\beta|}
\eeq
for every multi-indices $\alpha, \beta$. 
\begin{remark} We assume $|\xi|>1$ in the following. The case $|\xi|\leq1$ can be easily deduced from the regarding estimates.
\end{remark}
Let $\phi_i\in\mathcal{C}^\infty_o(\R^{\N_i}), i=1,2,\ldots,n$ such that 
\bel{phi}
\phi_i(z_i)~=~1,~~~~ |z_i|~\leq~{1\over 2}\qquad\hbox{and}\qquad \phi_i(z_i)~=~0,~~~~|z_i|~>~1.
\eeq
Define
 \bel{P}
\P(x,\xi)~=~\int_{\R^\N}   \p_x^\betaup \Omega(x,z) \prod_{i=1}^n e^{-2\pi\i z_i\cdot\xi_i}\left(\rho_i(\xi)\right)^{|\alphaup_i|} z_i^{\alpha_i}\prod_{i=1}^n\phi_i\left(\rho_i(\xi)z_i\right) dz.
\eeq
Observe that 
\bel{product bump function} 
z_i^{\alphaup_i} \phi_i\left(z_i\right)  \exp \Bigg\{ -2\pi\i  {z_i \cdot \xi_i \over \rho_i(\xi)}\Bigg\},\qquad i=1,2,\ldots,n
\eeq
is another normalized {\it bump}-function. 

By using (\ref{Cancellation}), we have
\bel{P norm}
\left|\P(x,\xi)\right|~\leq~\C_\beta~\left\{1+\sum_{i=1}^n ~\rho_i(\xi)\right\}^{\rho|\beta|}~\leq~\C_\beta~\left(1+|\xi|\right)^{\rho|\beta|}.
 \eeq
On the other hand,   define
\bel{E}
\E(x,\xi)~=~\int_{\R^\N}   \p_x^\betaup \Omega(x,z) \prod_{i=1}^n e^{-2\pi\i z_i\cdot\xi_i}\left(\rho_i(\xi)\right)^{|\alphaup_i|} z_i^{\alpha_i}\left\{1-\prod_{i=1}^n\phi_i\left(\rho_i(\xi)z_i\right)\right\} dz.
\eeq
Recall the multi-index $\gamma$ from the beginning of this section. 
Denote $|\xi_\ell|=\max_{i\in\{1,2,\ldots,n\}}|\xi_i|$.
Consider a series of integration by parts as follows.

We integrate by parts $w.r.t~z_\ell$ inside (\ref{E}) and stop  
whenever there is a partial derivative of $1-\prod_{i=1}^n\varphi_i\left(\rho_i(\xi)z_i\right)$ appeared. The resulting term can be expressed as
\bel{expression}
\begin{array}{lr}\ds
\E_\gamma(x,\xi)~=~
\Cc_\gamma\int_{\R^\N}|\xi_{\ell}|^{-\gamma_\ell}\prod_{i=1}^n \left(\rho_i(\xi)\right)^{|\alphaup_i|}  \Bigg\{\p_{z_\ell}^{\gamma_\ell} \p_x^\betaup \Omega(x,z) z_\ell^{\alpha_\ell}\Bigg\} 
\\\\ \ds~~~~~~~~~~~~~~~~~
\prod_{i=1}^n e^{-2\pi\i z_i\cdot\xi_i}\prod_{i\neq\ell} z_i^{\alpha_i}|\xi_\ell|^{-1}\p_{z_\ell}\left\{1-\prod_{i=1}^n\phi_i\left(z_i\right)\right\} dz
\end{array}
\eeq
where boundary terms vanish because $\p_{z_\ell}^{\gamma_\ell}\p_x^\beta\Omega(x,z)$ decays rapidly as $|z|\mt\infty$.

Note that $\prod_{i=1}^n\phi_i\left(\rho_i(\xi)z_i\right)=0$ if $|z_i|>\rho_i^{-1}(\xi)$ for any $i=1,2,\ldots,n$.
Therefore, we have $|\xi_\ell|^{-1}\p_{z_\ell}\left\{1-\prod_{i=1}^n\phi_i\left(\rho_i(\xi)z_i\right)\right\}$ as another product of normalized {\it bump}-functions
having the same characteristic of $\prod_{i=1}^n\phi_i(\rho_i(\xi)z_i)$. 

By carrying out a similar argument as (\ref{P})-(\ref{P norm}) and using (\ref{Cancellation}), we find
\bel{E gamma norm}
\left| \E_\gamma(x,\xi)\right|~\leq~\C_{\beta~\gamma}~|\xi|^{-\gamma_\ell} \left\{1+\sum_{i=1}^n ~\rho_i(\xi)\right\}^{\rho|\beta|+\gamma_\ell}~\leq~\C_{\beta~\gamma}~\left(1+|\xi|\right)^{\rho|\beta|}.
\eeq
We continue  this process until  $\gamma_\ell>\sum_{i=1}^n|\alpha_i|$.
Our assertion now reduces to
\bel{E rewrite}
\begin{array}{lr}\ds
\mathcal{E}(x,\xi)~=~\int_{\R^\N}|\xi_{\ell}|^{-\gamma_\ell}\prod_{i=1}^n \left(\rho_i(\xi)\right)^{|\alphaup_i|}  \Bigg\{\p_{z_\ell}^{\gamma_\ell} \p_x^\betaup \Omega(x,z) z_\ell^{\alpha_\ell}\Bigg\} 
\\\\ \ds~~~~~~~~~~~~~~~~~~~
\prod_{i=1}^n e^{-2\pi\i z_i\cdot\xi_i}\prod_{i\neq\ell} z_i^{\alpha_i}\left\{1-\prod_{i=1}^n\phi_i\left(\rho_i(\xi)z_i\right)\right\} dz
\\\\ \ds
~=~\int_{\R^\N}|\xi_{\ell}|^{-\gamma_\ell}\prod_{i=1}^n \left(\rho_i(\xi)\right)^{|\alphaup_i|}  \Bigg\{\p_{z_\ell}^{\gamma_\ell} \p_x^\betaup \Omega(x,z) z_\ell^{\alpha_\ell}\Bigg\} 
\\\\ \ds~~~~~~~~~~~~~~~~~~~~~~
\prod_{i=1}^n e^{-2\pi\i z_i\cdot\xi_i}\prod_{i\neq\ell} z_i^{\alpha_i}\left\{1-\prod_{i=1}^n\phi_i\left(\rho_i(\xi)z_i\right)\right\}\prod_{i=1}^n\phi(z_i) dz
\\\\ \ds
~+~\int_{\R^\N}|\xi_{\ell}|^{-\gamma_\ell}\prod_{i=1}^n \left(\rho_i(\xi)\right)^{|\alphaup_i|}  \Bigg\{\p_{z_\ell}^{\gamma_\ell} \p_x^\betaup \Omega(x,z) z_\ell^{\alpha_\ell}\Bigg\} 
\\\\ \ds~~~~~~~~~~~~~~~~~~~~~~
\prod_{i=1}^n e^{-2\pi\i z_i\cdot\xi_i}\prod_{i\neq\ell} z_i^{\alpha_i}\left\{1-\prod_{i=1}^n\phi_i\left(\rho_i(\xi)z_i\right)\right\}\left\{1-\prod_{i=1}^n\phi(z_i)\right\} dz
\\\\ \ds
~\doteq~\mathcal{E}_1(x,\xi)+\mathcal{E}_2(x,\xi).
\end{array}
\eeq
Note that $1-\prod_{i=1}^n\phi(z_i)=0$ if $|z_i|\leq{1\over 2},i=1,2,\ldots,n$. Hence that 
\bel{E_2 norm}
\left|\mathcal{E}_2(x,\xi)\right|~\leq~\C_{\alpha~\beta}
\eeq
 provided that $\p_{z_\ell}^{\gamma_\ell} \p_x^\betaup \Omega(x,z)$ decays rapidly as $|z|\mt\infty$.

Define
\bel{E_0}
\begin{array}{lr}\ds
\mathcal{E}_o(x,\xi)~=~\int_{\R^\N}|\xi_{\ell}|^{-\gamma_\ell}\prod_{i=1}^n \left(\rho_i(\xi)\right)^{|\alphaup_i|}  \Bigg\{\p_{z_\ell}^{\gamma_\ell} \p_x^\betaup \Omega(x,z) z_\ell^{\alpha_\ell}\Bigg\} 
\\\\ \ds~~~~~~~~~~~~~~~~~~~~~~
\prod_{i=1}^n e^{-2\pi\i z_i\cdot\xi_i}\prod_{i\neq\ell} z_i^{\alpha_i}\prod_{i=1}^n\phi_i\left(\rho_i(\xi)z_i\right)\prod_{i=1}^n\phi(z_i) dz
\\\\ \ds~~~~~~~~~~~~
~=~\int_{\R^\N}|\xi_{\ell}|^{-\gamma_\ell}\prod_{i=1}^n \left(\rho_i(\xi)\right)^{|\alphaup_i|}  \Bigg\{\p_{z_\ell}^{\gamma_\ell} \p_x^\betaup \Omega(x,z) z_\ell^{\alpha_\ell}\Bigg\} 
\\\\ \ds~~~~~~~~~~~~~~~~~~~~~~
\prod_{i=1}^n e^{-2\pi\i z_i\cdot\xi_i}\prod_{i\neq\ell} z_i^{\alpha_i}\prod_{i=1}^n\phi_i\left(\rho_i(\xi)z_i\right) dz.
\end{array}
\eeq
Again, by carrying out a similar argument as (\ref{P})-(\ref{P norm}) and using (\ref{Cancellation}), we find
\bel{E 0 norm}
\left| \mathcal{E}_o(x,\xi)\right|~\leq~\C_{\beta~\gamma}~|\xi|^{-\gamma_\ell} \left\{1+\sum_{i=1}^n ~\rho_i(\xi)\right\}^{\rho|\beta|+\gamma_\ell}~\leq~\C_{\alpha~\beta}~\left(1+|\xi|\right)^{\rho|\beta|}.
\eeq
On the other hand, the norm of
\bel{E_3}
\begin{array}{lr}\ds
\mathcal{E}_3(x,\xi)~=~\int_{\R^\N}|\xi_{\ell}|^{-\gamma_\ell}\prod_{i=1}^n \left(\rho_i(\xi)\right)^{|\alphaup_i|}  \Bigg\{\p_{z_\ell}^{\gamma_\ell} \p_x^\betaup \Omega(x,z) z_\ell^{\alpha_\ell}\Bigg\} 
\\\\ \ds~~~~~~~~~~~~~~~~~~~~~~
\prod_{i=1}^n e^{-2\pi\i z_i\cdot\xi_i}\prod_{i\neq\ell} z_i^{\alpha_i}\prod_{i=1}^n\phi_i\left(z_i\right) dz
\end{array}
\eeq
is bounded by a constant $\C_{\alpha~\beta}$.  Together with (\ref{E 0 norm}), we conclude that 
\bel{E 1 norm}
\left| \mathcal{E}_1(x,\xi)\right|~\leq~\C_{\alpha~\beta}~\left(1+|\xi|\right)^{\rho|\beta|}.
\eeq

\section{A fundamental lemma}
\setcounter{equation}{0}
Let $0<\rhoup<1$. Consider $\Theta\in\mathcal{C}^\infty(\R^\N\times\R^\N\times\R^\N\times\R^\N)$ of which
\bel{INEQ}
\begin{array}{lr}\ds
\left|\p_\xi^{\alphaup^1}\p_\eta^{\alphaup^2}\p_x^{\beta^1}\p_y^{\beta^2}\Theta(x,y,\xi,\eta)\right|
\\\\ \ds 
\leq~\C_{\alpha^1~\alpha^2~\beta^1~\beta^2}~\vartheta(\xi,\eta)
\prod_{i=1}^n\left({1\over 1+|\xi_i|+|\xi|^{\rhoup}}\right)^{|\alphaup_i^1|} \left({1\over 1+|\eta_i|+|\eta|^{\rhoup}}\right)^{|\alpha_i^2|}(1+|\xi|)^{\rhoup|\beta^1|}(1+|\eta|)^{\rhoup|\beta^2|}
\end{array}
\eeq
for every multi-indices $\alphaup^1,\alpha^2,\beta^1,\beta^2$.

Define
\bel{Lambda iint}
\Lambda(x,\xi)~=~\iint_{\R^\N\times\R^\N} e^{2\pi\i(x-y)\cdot(\xi-\eta)} \Theta(x,y,\xi,\eta) dy d\eta.
\eeq

{\bf Lemma One}~~{\it Let $\Lambda$ defined in (\ref{Lambda iint}).
Suppose that $\Theta\in\mathcal{C}^\infty(\R^\N\times\R^\N\times\R^\N\times\R^\N)$ 
satisfies the differential inequality in (\ref{INEQ}).
 We have
\bel{ineq}
\left|\p_\xi^\alphaup\p_x^\betaup\Lambda(x,\xi)\right|
~\leq~\C_{\alphaup~\betaup~\rhoup}~
\prod_{i=1}^n\left({1\over 1+|\xi_i|+|\xi|^{\rhoup}}\right)^{|\alphaup_i|}(1+|\xi|)^{\rhoup|\betaup|}
\eeq
for every multi-indices $\alphaup$, $\betaup$.}

\begin{cor} 
Denote $\sigma_i$ to be the symbol of $\T_i$ for $i=1,2$. Moreover, $\sigma$ is the symbol of
 $\T_1\circ \T_2$.  Suppose $\sigma_1,\sigma_2\in\S_\rhoup$. 
We have $\sigma\in\S_\rhoup$.
\end{cor}
{\bf Proof}~  A direct computation shows
\bel{symbol compo}
\sigma(x,\xi)~=~\iint_{\R^\N\times\R^\N} e^{2\pi\i(x-y)\cdot(\xi-\eta)}\sigma_1(x,\eta)\sigma_2(y,\xi)dyd\eta.
\eeq
Note that $\sigma_1(x,\eta)\sigma_2(y,\xi)$ satisfies the differential inequality in (\ref{INEQ}). 
By using {\bf Lemma One}, we have $\sigma$ in (\ref{symbol compo}) satisfying the differential inequality in (\ref{ineq}). 
\endproof

We develop the proof of {\bf Lemma One} in analogue to the work of Boutet de Monvel \cite{Monvel} and Beals and Fefferman \cite{Beals-Fefferman}.

Let $0<\rhoup<1$. Define a differential operator 
\bel{D}
\D~=~I-\left({1\over 4\pi^2}\right)\left({1\over 1+|\xi|^{2}}+{1\over 1+|\eta|^{2}}\right)^\rhoup\Delta_y-\left({1\over 4\pi^2}\right)\left(  1+|\xi|^2+|\eta|^2\right)^{\rhoup}\Delta_\eta
\eeq
and the regarding quadratic function
\bel{Q}
 \Q(x,y,\xi,\eta)~=~1+\left({1\over 1+|\xi|^{2}}+{1\over 1+|\eta|^{2}}\right)^\rhoup|\xi-\eta|^{2}+\left( 1+|\xi|^2+ |\eta|^2\right)^{\rhoup}|x-y|^{2}.
 \eeq
Denote $\L=\Q^{-1}\D$. We find 
\bel{L identity}
\L^N \exp\Bigg\{2\pi\i(x-y)\cdot(\xi-\eta)\Bigg\}~=~\exp\Bigg\{2\pi\i(x-y)\cdot(\xi-\eta)\Bigg\},\qquad N\ge1.
\eeq
\begin{remark}
We momentarily assume that $\Theta(x,y,\xi,\eta)$ has compact support in both $y$ and $\eta$.
\end{remark}
By using the identity in (\ref{L identity}) and integration 
 by parts $w.r.t$ $y$ and $\eta$ inside (\ref{Lambda iint}), we obtain
\bel{Integral'}
\iint_{\R^\N\times\R^\N} {^t}\L^N\Theta(x,y,\xi,\eta)e^{2\pi\i(x-y)\cdot(\xi-\eta)}dyd\eta
\eeq
where ${^t}\L={^t}\D\Q^{-1}$. 
\begin{lemma} Suppose that $\Theta\in\mathcal{C}^\infty(\R^\N\times\R^\N\times\R^\N\times\R^\N)$ satisfies (\ref{INEQ}). We have
\bel{phi est}
\left|{^t}\L^N\Theta(x,y,\xi,\eta)\right|~\leq~\C_N~\vartheta(\xi,\eta) \left(1+{|\xi|^2\over 1+|\eta|^2}+{|\eta|^2\over 1+|\xi|^2}\right)^{\rhoup N} \Q^{-N}(x,y,\xi,\eta)
\eeq
 for every $N\ge1$.
\end{lemma}
{\bf Proof }~ Denote 
\bel{a,b}
\a(\xi,\eta)~=~\left({1\over 1+|\xi|^{2}}+{1\over 1+|\eta|^{2}}\right)^\rhoup
,\qquad \b(\xi,\eta)~=~ \left(1+|\xi|^2+ |\eta|^2\right)^{\rhoup}.
\eeq
From (\ref{D})-(\ref{Q}), we find
\bel{form}
{^t}\L^N \Theta(x,y,\xi,\eta)~=~\left( {1\over \Q} -\Delta_y {\a\over \Q} -\Delta_\eta {\b\over \Q}\right)^N \Theta(x,y,\xi,\eta),\qquad N\ge1.
\eeq
Observe that
\bel{es ab}\begin{array}{cc}\ds
\a(\xi,\eta)\Big|\Delta_y \Theta(x,y,\xi,\eta)\Big|~\leq~\C~\a(\xi,\eta)\left(1+|\eta|\right)^{2\rhoup} \vartheta(\xi,\eta)
\\\\ \ds~~~~~~~~~~~~~~~~~~~~~~~~~~~~~~~~~~~~~~~~~~~~~~~~~~~
~\leq~\C~\left(1+{|\xi|^2\over 1+|\eta|^2}+{|\eta|^2\over 1+|\xi|^2}\right)^{\rhoup}\vartheta(\xi,\eta),
\\\\ \ds
\b(\xi,\eta)\Big|\Delta_\eta \Theta(x,y,\xi,\eta)\Big|~\leq~\C~\b(\xi,\eta) \left({1\over 1+|\eta|}\right)^{2\rhoup}\vartheta(\xi,\eta)
\\\\ \ds~~~~~~~~~~~~~~~~~~~~~~~~~~~~~~~~~~~~~~~~~~~~~~~~~~~
~\leq~\C~\left(1+{|\xi|^2\over 1+|\eta|^2}+{|\eta|^2\over 1+|\xi|^2}\right)^{\rhoup}\vartheta(\xi,\eta).
\end{array}
\eeq
Hence that if all partial derivatives fall on $\Theta$ in (\ref{form}), then (\ref{es ab}) implies (\ref{phi est}).

Turn to the general case.   From (\ref{Q}) and (\ref{a,b}), we verify  that
\bel{est2.3}
\left|\p_\eta^\alpha \a(\xi,\eta)\right|~\leq~\C_\alpha~\left({1\over 1+|\eta|}\right)^{|\alphaup|}\a(\xi,\eta)
,\qquad
\left|\p_\eta^\alphaup \b(\xi,\eta)\right|~\leq~\C_\alpha~\left({1\over 1+|\xi|+|\eta|}\right)^{|\alphaup|}\b(\xi,\eta)
\eeq
for every multi-index $\alphaup$. 

Our task can be finished if we show
\bel{est Crucial} 
\left| \p_\eta^\alphaup\p_y^\betaup \Q^{-1}(x,y,\xi,\eta)\right|~\leq~\C_{\alpha~\beta}~ \Q^{-1}(x,y,\xi,\eta)\left({1\over 1+|\eta|}\right)^{\rhoup|\alphaup|}\left(1+|\xi|+|\eta|\right)^{\rhoup|\betaup|}
\eeq
for every multi-indices $\alphaup$, $\betaup$.  
 
Note that  $\p_\eta^\alphaup\p_y^\betaup \Q^{-1}$ consists a linear combination of
\bel{Expansion}
\Q^{-1}\prod_k \left[\Q^{-1}\p_\eta^{\alphaup^k}\p_y^{\betaup^k}\Q\right]^k,\qquad  |\alphaup|~=~\sum_k |\alphaup^k|k,\qquad |\betaup|=\sum_k|\betaup^k|k.
\eeq
We need the following preliminary estimates:

$\bullet$ Suppose $|\xi-\eta|\leq(1+|\xi|+|\eta|)^\rhoup$. We necessarily have $|\xi|\sim|\eta|$ and 
\bel{pre est1}
\begin{array}{cc}\ds
 \left({1\over 1+|\xi|^2}+{1\over 1+|\eta|^2}\right)^{\rhoup}|\xi-\eta|
~\leq~\C~\left({1\over 1+|\eta|}\right)^\rhoup,
\\\\ \ds
 \left({1\over 1+|\xi|^2}+{1\over 1+|\eta|^2}\right)^{\rhoup}
~\leq~\C~\left({1\over 1+|\eta|}\right)^{2\rhoup}.
\end{array}
\eeq
 $\bullet$ Suppose $|\xi-\eta|>(1+|\xi|+|\eta|)^\rhoup$. We have
 \bel{pre est2}
 \begin{array}{cc}\ds
 \Q^{-1}(x,y,\xi,\eta) \left({1\over 1+|\xi|^2}+{1\over 1+|\eta|^2}\right)^{\rhoup}|\xi-\eta|
 \\\\ \ds
 ~\leq~\C |\xi-\eta|^{-1}
~\leq~\C \left({1\over 1+|\eta|}\right)^\rhoup,
\\\\ \ds
\Q^{-1}(x,y,\xi,\eta) \left({1\over 1+|\xi|^2}+{1\over 1+|\eta|^2}\right)^{\rhoup}
~\leq~\C |\xi-\eta|^{-2}
~\leq~\C \left({1\over 1+|\eta|}\right)^{2\rhoup}.
\end{array}
\eeq
$\bullet$ Suppose $|x-y|\leq(1+|\xi|+|\eta|)^{-\rhoup}$.  We have
\bel{pre est3}
\begin{array}{cc}\ds
 \left(1+|\xi|^2+|\eta|^2\right)^\rhoup|x-y|~\leq~\C\left(1+|\xi|+|\eta|\right)^\rhoup,
 \\\\ \ds
  \left(1+|\xi|^2+|\eta|^2\right)^\rhoup~\leq~\C\left(1+|\xi|+|\eta|\right)^{2\rhoup}.
  \end{array}
  \eeq
$\bullet$ Suppose $|x-y|>(1+|\xi|+|\eta|)^{-\rhoup}$.  We have
\bel{pre est4}
\begin{array}{cc}\ds
\Q^{-1}(x,y,\xi,\eta) \left(1+|\xi|^2+|\eta|^2\right)^\rhoup|x-y|
~\leq~\C|x-y|^{-1}
~\leq~\C\left(1+|\xi|+|\eta|\right)^\rhoup,
\\\\ \ds
\Q^{-1}(x,y,\xi,\eta) \left(1+|\xi|^2+|\eta|^2\right)^\rhoup
~\leq~\C|x-y|^{-2}
~\leq~\C\left(1+|\xi|+|\eta|\right)^{2\rhoup}.
\end{array}
\eeq

From direct computation and by using  (\ref{pre est1})-(\ref{pre est4}) together with (\ref{est2.3}), we find 
\bel{Building block Result}
\left|\Q^{-1}\p_\eta^\alphaup\p_y^\betaup \Q(x,y,\xi,\eta)\right|~\leq~\C_{\alpha~\beta}~\left({1\over 1+|\eta|}\right)^{\rhoup|\alphaup|}\left(1+|\xi|+|\eta|\right)^{\rhoup|\betaup|}
\eeq
for every multi-indices $\alphaup$, $\betaup$.
\endproof

{\bf Proof of Lemma One} ~Recall (\ref{Lambda iint}). 
$\p_\xi^\alphaup\p_x^\betaup\Lambda(x,\xi)$ can be expressed as
 a linear combination of  
\bel{diff term}
\begin{array}{lr}\ds
\iint_{\R^\N\times\R^\N} e^{2\pi\i(x-y)\cdot(\xi-\eta)} \p_\xi^{\alphaup^1} \p_x^{\betaup^1}    \p_\eta^{\alphaup^2} \p_y^{\betaup^2}  \Theta(x,y,\xi,\eta) dy d\eta
\end{array}
\eeq
where $|\alphaup_i|=|\alphaup^1_i|+|\alphaup^2_i|, i=1,2,\ldots,n$ and $|\betaup|=|\betaup^1|+|\betaup^2|$. 

Recall (\ref{L identity}) and ${^t}\L={^t}\D\Q^{-1}$ as shown in (\ref{Integral'}). By integration by parts $w.r.t$ $y$ and $\eta$, we find  (\ref{diff term}) equal to
\bel{diff term rewrite}
\iint_{\R^\N\times\R^\N} e^{2\pi\i(x-y)\cdot(\xi-\eta)}{^t}\L^N \p_\xi^{\alphaup^1} \p_x^{\betaup^1}    \p_\eta^{\alphaup^2} \p_y^{\betaup^2}  \Theta(x,y,\xi,\eta) dy d\eta,\qquad N\ge1.
\eeq
{\bf Case One}~~Suppose $|\xi-\eta|>{1\over 2}(1+|\xi|+|\eta|)$. We carry out an $M$-fold integration by parts $w.r.t~y$ inside (\ref{diff term rewrite}). 

The  function $\Delta_y^M \p_\xi^{\alphaup^1}\p_\eta^{\alphaup^2}\p_x^{\betaup^1}\p_y^{\betaup^2}\Theta(x,y,\xi,\eta)/|\xi-\eta|^{2M}$
satisfies the differential inequality in (\ref{INEQ}) with norm bounded by 
\bel{Delta d varphi norm}
\begin{array}{lr}\ds
\C_{\alpha~\beta~M} ~\left({1+|\xi|+|\eta|}\right)^{2(\rhoup-1)M}
\\ \ds
\prod_{i=1}^n\left({1\over 1+|\xi_i|+|\xi|^{\rhoup}}\right)^{|\alphaup_i^1|} \left({1\over 1+|\eta_i|+|\eta|^{\rhoup}}\right)^{|\alphaup_i^2|}(1+|\xi|)^{\rhoup|\betaup^1|}(1+|\eta|)^{\rhoup|\betaup^2|}.
\end{array}
\eeq
By applying {\bf Lemma 4.1}, we find
\bel{L^t Delta d varphi norm}
\begin{array}{lr}\ds
\left|^t\L^N\left({1\over |\xi-\eta|}\right)^{2M}\Delta_y^M \p_\xi^{\alphaup^1}\p_\eta^{\alphaup^2}\p_x^{\betaup^1}\p_y^{\betaup^2}\Theta(x,y,\xi,\eta)\right|
\\\\ \ds
~\leq~\C_{\alpha~\beta~M~N}~ \left(1+{|\xi|^2\over 1+|\eta|^2}+{|\eta|^2\over 1+|\xi|^2}\right)^{\rhoup N} \Q^{-N}(x,y,\xi,\eta)
\\\\ \ds
~~~\left({1+|\xi|+|\eta|}\right)^{2(\rhoup-1)M}\prod_{i=1}^n\left({1\over 1+|\xi_i|+|\xi|^{\rhoup}}\right)^{|\alphaup_i^1|} \left({1\over 1+|\eta_i|+|\eta|^{\rhoup}}\right)^{|\alphaup_i^2|}(1+|\xi|)^{\rhoup|\betaup^1|}(1+|\eta|)^{\rhoup|\betaup^2|}.
\end{array}
\eeq
By choosing $M$ sufficiently large, depending on $\alphaup$, $\betaup$ and  $\rhoup$,  we have
\bel{diff term norm2}
\begin{array}{lr}\ds
 \iint_{|\xi-\eta|~>~{1\over 2}(1+|\xi|+|\eta|)} e^{2\pi\i(x-y)\cdot(\xi-\eta)} {^t}\L^N \left({1\over |\xi-\eta|}\right)^{2M}\Delta^M_y \p_\xi^{\alphaup^1} \p_x^{\betaup^1}    \p_\eta^{\alphaup^2} \p_y^{\betaup^2}  \Theta(x,y,\xi,\eta) dy d\eta
\\\\ \ds
~\leq~\C_{\alpha~\betaup~\rhoup~ N}~ \prod_{i=1}^n\left({1\over 1+|\xi_i|+|\xi|^{\rhoup}}\right)^{|\alphaup_i|}(1+|\xi|)^{\rhoup|\betaup|}\iint_{\R^\N\times\R^\N} \Q^{-N}(x,y,\xi,\eta)dyd\eta.
\end{array}
\eeq
Recall $\Q(x,y,\xi,\eta)$ defined in (\ref{Q}).  By changing variables 
$\eta \mt \left(1+|\xi|^2\right)^{\rhoup/2}(\xi-\eta)$ and $y\mt\left({1\over1+|\xi|^2}\right)^{\rhoup/2}(x-y)$, 
we find
\bel{est2.12}
\iint_{\R^\N\times\R^\N} \Q^{-N}(x,y,\xi,\eta)dyd\eta~\leq~\C~\iint_{\R^\N\times\R^\N}  \Big(1+|\eta|^{2}+|y|^{2}\Big)^{-N} dyd\eta.
\eeq
The integral converges if $N$ is sufficiently large.

{\bf Case Two} ~~Suppose $|\xi-\eta|\leq{1\over 2}(1+|\xi|+|\eta|)$. We necessarily have  $|\xi|\sim|\eta|$.
Let $\U\cup\V=\{1,2,\ldots,n\}$ such that
$|\xi_i-\eta_i|\leq{1\over 2}(1+|\xi_i|+|\eta_i|)$ if $i\in\U$ and $|\xi_i-\eta_i|>{1\over 2}(1+|\xi_i|+|\eta_i|)$ if $i\in\V$.

 For every $i\in\U$,  we have  $|\xi_i|\sim|\eta_i|$. Suppose $|\xi_i|\lesssim|\xi|^\rho$ for every $i\in\V$. 
 The function $\p_\xi^{\alphaup^1}\p_\eta^{\alphaup^2}\p_x^{\betaup^1}\p_y^{\betaup^2}\Theta(x,y,\xi,\eta)$ satisfies the differential inequality in (\ref{INEQ}) with  norm bounded by 
 \bel{d Theta norm}
\begin{array}{lr}\ds
\C_{\alphaup~\betaup}~\prod_{i\in\U}\left({1\over 1+|\xi_i|+|\xi|^{\rhoup}}\right)^{|\alphaup_i|}\prod_{i\in\V}\left({1\over 1+|\xi_i|+|\xi|^{\rhoup}}\right)^{|\alphaup_i^1|}\left({1\over 1+|\eta_i|+|\eta|^{\rhoup}}\right)^{|\alphaup_i^2|}(1+|\xi|)^{\rhoup|\betaup|}
\\\\ \ds
~\leq~\C_{\alphaup~\betaup}~\prod_{i=1}^n\left({1\over 1+|\xi_i|+|\xi|^{\rhoup}}\right)^{|\alphaup_i|}(1+|\xi|)^{\rhoup|\betaup|}.
\end{array}
\eeq
 On the other hand, assume that there is at least one $i\in\V$ such that $|\xi|^\rho\lesssim|\xi_i|$. 
 Denote $y_\vv$ to be the projection of $y$ on the subspace $\bigotimes_{i\in\V}\R^{\N_i}$, same for $\xi_\vv$ and $\eta_\vv$.
 
Again,  we carry out an $M$-fold integration by parts $w.r.t~y$ inside (\ref{diff term rewrite}). 
The function $\Delta_{y_\vv}^M \p_\xi^{\alphaup^1}\p_\eta^{\alphaup^2}\p_x^{\betaup^1}\p_y^{\betaup^2}\Theta(x,y,\xi,\eta)/|\xi_\vv-\eta_\vv|^{2M}$ satisfies the differential inequality in (\ref{INEQ}) with  norm bounded by
\bel{d Theta Norm}
\begin{array}{lr}\ds
\C_{\alphaup~\betaup~M }~\left({1+|\xi_\vv|+|\eta_\vv|}\right)^{-2M}\left(1+|\xi|^\rho\right)^{2M}
\\\\ \ds
\prod_{i\in\U}\left({1\over 1+|\xi_i|+|\xi|^{\rhoup}}\right)^{|\alphaup_i|}\prod_{i\in\V}\left({1\over 1+|\xi_i|+|\xi|^{\rhoup}}\right)^{|\alphaup_i^1|}\left({1\over 1+|\eta_i|+|\eta|^{\rhoup}}\right)^{|\alphaup_i^2|}(1+|\xi|)^{\rhoup|\betaup|}
\\\\ \ds
~\leq~\C_{\alphaup~\betaup~M }~\left({1+|\xi_\vv|+|\eta_\vv|}\right)^{-2M}\left(1+|\xi|^\rho\right)^{2M}
\\\\ \ds
\prod_{i\in\U}\left({1\over 1+|\xi_i|+|\xi|^{\rhoup}}\right)^{|\alphaup_i|}\prod_{i\in\V}\left({1\over 1+|\xi_i|+|\xi|^{\rhoup}}\right)^{|\alphaup_i^1|}\left({1\over 1+|\xi|^{\rhoup}}\right)^{|\alphaup_i^2|}(1+|\xi|)^{\rhoup|\betaup|}
\\\\ \ds
~\leq~\C_{\alpha~\beta}~\prod_{i=1}^n\left({1\over 1+|\xi_i|+|\xi|^{\rhoup}}\right)^{|\alphaup_i|} (1+|\xi|)^{\rhoup|\betaup|}
\end{array}
\eeq
provided that $M$ is chosen sufficiently large depending on $\alpha$.

By repeating the estimates in analogue to (\ref{L^t Delta d varphi norm})-(\ref{est2.12}), we conclude the same result for $|\xi-\eta|\leq{1\over 2}(1+|\xi|+|\eta|)$.

Our estimates above are independent from the size of $\supp\Theta$ in $y$ and $\eta$. The assumption of compactness can be removed by taking the approximation as discussed in {\bf 1.3}, chapterVI of Stein \cite{S}. 
\endproof

\section{Point-wise estimates on partial operators}
\setcounter{equation}{0}
Denote $\M$ to be the strong maximal operator. Recall  $\Delta_\t$ defined in (\ref{Delta_t}). 
\begin{lemma}  Let $\sigma\in\S_\rhoup, 0<\rho<1$. We have
\bel{Delta_t Tf Est}
\left|\Delta_\t \T f(x)\right|~\leq~\C_\rho~\M f(x)
\eeq
for every $\t\in\H$  satisfying (\ref{H}).
\end{lemma}

{\bf Proof :} A direct computation shows 
\bel{Delta_t T}
\begin{array}{lr}\ds
\Delta_\t \T f(x)
~=~\int_{\R^\N} f(y)\Omega_\t(x,y)dy
\\\\ \ds ~~~~~~~~~~~~~
~=~\int_{\R^\N} f(y) \left\{  \int_{\R^\N} e^{2\pi\i (x-y)\cdot\xi}\deltaup_\t(\xi) \Lambda(y,\xi)d\xi\right\}dy
\end{array}
\eeq
where
\bel{Lambda}
\Lambda(x,\xi)~=~\iint_{\R^\N\times\R^\N} e^{2\pi\i(x-y)\cdot(\xi-\eta)}\sigma(y,\eta)dyd\eta.
\eeq
Let $\sigma\in\S_\rhoup, 0<\rho<1$. 
By applying {\bf Lemma One}, we find that $\Lambda(x,\xi)$ is bounded and satisfies the differential inequality in (\ref{ineq}). 

Denote $z=x-y$ and 
 \bel{t dila}
\begin{array}{lr}\ds
\t\xi~=~\left(2^{-t_1}\xi_1,~2^{-t_2}\xi_2,\ldots,~2^{-t_n}\xi_n\right),
\qquad
\t^{-1}\xi~=~\left(2^{t_1}\xi_1,~2^{t_2}\xi_2,\ldots,~2^{t_n}\xi_n\right).
\end{array}
\eeq
By changing dilations $\xi\mt\t^{-1}\xi$ and $z\mt\t z$, the integral in (\ref{Delta_t T}) 
can be written as 
\bel{Delta_k T dila}
\int_{\R^\N} f(x-\t z) \left\{ \int_{\R^\N} e^{2\pi\i z\cdot\xi} \deltaup_\t\left(\t^{-1}\xi\right) \Lambda( x-\t z,\t^{-1}\xi)d\xi\right\}dz.
\eeq
Let $q=1/\rhoup$. By definition of $\deltaup_\t(\xi)$ in (\ref{delta_t}), the support of 
\bel{delta dila-change}
\deltaup_\t\left(\t^{-1}\xi\right)~=~\prod_{i=1}^n \phi\left(2^{t_1-q t_i}\xi_1,\ldots,\xi_i,\ldots,2^{t_n-q t_i}\xi_n\right)
\eeq
is contained  in a ball with radius equal to $2$.

Let $\I\cup\J=\{1,2,\ldots,n\}$ defined in (\ref{IJ}).  By integration by parts $w.r.t~\xi$, we find
\bel{Delta_k T dila'}
\begin{array}{lr}\ds
 \int_{\R^\N} e^{2\pi\i z\cdot\xi} \deltaup_\t\left(\t^{-1}\xi\right) \Lambda(x-\t z,\t^{-1}\xi)d\xi
\\\\ \ds
~=~\left({1\over 2\pi\i z}\right)^{|\alpha|}  \int e^{2\pi\i z\cdot\xi} \prod_{i\in\I}\p_{\xi_i}^{\alphaup_i}\prod_{i\in\J}\p_{\xi_i}^{\alphaup_i}\deltaup_\t\left(\t^{-1}\xi\right) \Lambda(x-\t z,\t^{-1}\xi)d\xi
\end{array}
\eeq
at $z_i\neq0$, $i=1,2,\ldots,n$ for every multi-index $\alphaup$.

The  integral in (\ref{Delta_k T dila'}) has a  bounded norm because $\left|\supp\deltaup_\t\left(\t^{-1}\xi\right)\right|\leq\C 2^n$. 

Recall {\bf Lemma 2.1}. 
Given $\t\in\H$ satisfying (\ref{H}), we have $|\xi_i|\sim 1$ for  $i\in\I$ whenever $\xi\in\supp\deltaup_\t\left(\t^{-1}\xi\right)$. In particular,   $\imath\in\I$ of which $|\xi_\imath|=\max_{i\in\{1,2,\ldots,n\}}|\xi_i|$ and $2^{t_\imath}\sim2^{q t_i}$ for $i\in\J$.

From direct computation, we find
\bel{Est Case}
\begin{array}{lr}\ds
\left| \prod_{i\in\I}\p_{\xi_i}^{\alphaup_i}\prod_{i\in\J}\p_{\xi_i}^{\alphaup_i} \deltaup_\t\left(\t^{-1}\xi\right)  \Lambda(x-\t z,\t^{-1}\xi)\right|
\\\\ \ds 
~\leq~\C_\rho\prod_{i\in\I}\left(2^{t_i}\right)^{|\alphaup_i|} \left({1\over 1+2^{t_i}|\xi_i|+|\t^{-1}\xi|^{\rhoup}}\right)^{|\alphaup_i|} \prod_{i\in\J} \left(2^{t_i}\right)^{|\alphaup_i|}\left({1\over 1+2^{t_i}|\xi_i|+|\t^{-1}\xi|^{\rhoup}}\right)^{|\alphaup_i|}
\\\\ \ds 
~\leq~\C_\rho\prod_{i\in\I}\left(2^{t_i}\right)^{|\alphaup_i|} \left({1\over 1+2^{t_i}|\xi_i|}\right)^{|\alphaup_i|} \prod_{i\in\J} \left(2^{t_i}\right)^{|\alphaup_i|}\left({1\over 1+2^{t_\imath\rho}|\xi_\imath|^{\rhoup}}\right)^{|\alphaup_i|}
\\\\ \ds 
~\leq~\C_\rho\prod_{i\in\I}\left({1\over |\xi_i|}\right)^{|\alphaup_i|} \prod_{i\in\J} \left({1\over |\xi_\imath|}\right)^{\rhoup|\alphaup_i|}
\\\\ \ds 
~\leq~\C_{\alpha~\rho}.
\end{array}
\eeq
From (\ref{Delta_k T dila'})-(\ref{Est Case}), we obtain
\bel{Delta_t Tf Est Int}
\left|\Delta_\t \T f(x)\right|~\leq~ \C_{\rho~N}\int |f(x-\t z)|\left({1\over 1+|z|}\right)^{N} dz,\qquad N\ge1.
\eeq
Note that $\Big( 1+|z|\Big)^{-N}$ can be approximated by $\sum_{k=0}^\infty a_k\chi_{\B_k}$ where  $a_k>0$  and $\chi_{\B_k}$  is the indicator function on a ball $\B_k$ centered on the origin of $\R^{\N}$ with radius $r_k>0$.  Moreover, 
$\sum_k a_k|\B_k|<\infty$ provided that  $N$ is sufficiently large. 
From (\ref{Delta_t Tf Est Int}), we have
\bel{Delta_t f Est Max}
\begin{array}{lr}\ds
\left|\Delta_\t \T f(x)\right|~\leq~ \C_\rho\int_{\R^\N} |f(x-\t z)|\sum_{k=0}^\infty a_k\chi_{\B_k} dz
\\\\ \ds~~~~~~~~~~~~~~~~
~\leq~\C_\rho\left\{\sum_{k=0}^\infty a_k|\B_k|\right\} \sup_k {1\over|\B_k|}\int_{\B_k} |f(x-\t z)|dz
\\\\ \ds~~~~~~~~~~~~~~~~
~\leq~\C_\rho~\M f(x).
\end{array}
\eeq
\endproof

Denote $\s$ to be another $n$-tuple  as same as $\t$.

\begin{lemma}   Suppose $\sigma\in\S_\rhoup, 0<\rho<1$. We have
\bel{decay Est}
\left| \Delta_\t \T \Delta_\s f(x)\right|~\leq~\C_\rho
\prod_{i=1}^n2^{-(1-\rhoup)|t_i-s_i|}~\M f(x)
\eeq
for every $\t, \s\in\H$ satisfying (\ref{H}).
\end{lemma}

{\bf Proof} ~~Recall (\ref{Delta_t T})-(\ref{Lambda}). 
From direct computation,  we have
\bel{Delta_t T Delta_s f}
\begin{array}{lr}\ds
 \Delta_\t\T \Delta_\s f(x)~=~\int_{\R^\N} f(y)\Omega_{\t\s}(x,y)dy
\\\\ \ds ~~~~~~~~~~~~~~~~~
~=~\int_{\R^\N} f(y)\left\{\int_{\R^\N} e^{2\pi\i (x-y)\cdot\xi} \deltaup_\t(\xi)\Lambda_\s(y,\xi)d\xi\right\}dy
\end{array}
\eeq
where
\bel{Lambda_s}
\Lambda_\s(x,\xi)~=~
\iint_{\R^\N\times\R^\N} e^{2\pi\i(x-y)\cdot(\xi-\eta)} \deltaup_\s(\eta)\sigma(y,\eta)dyd\eta.
\eeq
In the following paragraph, we write $\c>0$ for some fixed, suitable  constant.

Define $|\t-\s|=\max_{i\in\{1,2,\ldots,n\}} |t_i-s_i|$. Recall {\bf Lemma 2.2}. Note that  $\deltaup_\s(\eta)\sigma(y,\eta)$ satisfies the differential inequality in (\ref{INEQ}).
By using {\bf Lemma One}, we find that $\Lambda_\s(x,\xi)$ is bounded and satisfies the differential inequality in (\ref{ineq}). For $|\t-\s|\leq\c$, the situation can be handled as {\bf Lemma 5.1}.

From now on, we consider $|\t-\s|>\c$.
There is at least one $i\in\{1,2,\ldots,n\}$ such that $|t_i-s_i|>\c$.
By integration by parts $w.r.t~y_\ell$, we find
\bel{int part Lambda_s}
\Lambda_\s(x,\xi)~=~{-1\over 4\pi^2}
\iint_{\R^\N\times\R^\N} e^{2\pi\i(x-y)\cdot(\xi-\eta)}\deltaup_\s(\eta) \Delta_{y_i}\sigma(y,\eta){1\over |\xi_i-\eta_i|^2}dyd\eta.
\eeq
Note that the boundary terms are vanished. This can be seen by integration by parts $w.r.t~\eta$ whereas $\delta_\s(\eta) \Delta_{y_i}\sigma(y,\eta)/ |\xi_i-\eta_i|^2$ has a compact support in $\eta$ for every $\s$.

Let $\sigma\in\S_\rhoup, 0<\rho<1$ satisfying (\ref{Ineq}) and $q=1/\rho$. Given $\t,\s\in\H$ as (\ref{H}),  we consider $\xi\in\supp \deltaup_\t(\xi)$ and $\eta\in\supp \deltaup_\s(\eta)$. 

Denote $t_\imath=\max_{i\in\{1,2,\ldots,n\}} t_i$ and $s_\jmath=\max_{i\in\{1,2,\ldots,n\}}s_i$.
We define
$\I_1\cup\J_1=\I_2\cup\J_2=\{1,2,\ldots,n\}$
  respectively for $\t$ and $\s$ as  (\ref{IJ}).

Recall {\bf Lemma 2.1}.  We have $|\xi_i|\sim2^{t_i}, i\in\I_1$ and $|\xi_i|\lesssim2^{t_i}, 2^{t_\imath}\sim2^{q t_i}, i\in\J_1$. On the other hand, we have 
$|\eta_i|\sim2^{s_i}, i\in\I_2$ and $|\eta_i|\lesssim2^{s_i}, 2^{s_\jmath}\sim2^{q s_i}, i\in\J_2$.
\begin{remark}   We must have $2^{t_i}\lesssim 2^{q t_j}$ and $2^{s_i}\lesssim2^{q s_j}$ for every $i,j\in\{1,2,\ldots,n\}$ provided that $\supp\delta_\t(\xi)\neq\emptyset$ and $\supp\delta_\s(\eta)\neq\emptyset$.
\end{remark}
{\bf 1.} Suppose $i\in\I_1\cap\I_2$. We have
$|\xi_i|\sim 2^{t_i}$,  $|\eta_i|\sim2^{s_i}$ and therefore
\bel{norm 2 cases}
|\xi_i-\eta_i|~\sim~\left\{\begin{array}{lr}\ds 2^{t_i}~=~2^{s_i}2^{(t_i-s_i)}, \qquad t_i>s_i+\c,
\\\\ \ds
2^{s_i}~=~2^{t_i}2^{(s_i-t_i)},\qquad t_i<s_i-\c.
\end{array}\right.
\eeq
 From {\bf Remark 5.1}, the first equality in (\ref{norm 2 cases}) further implies
\bel{iint norm est2}
{1+|\eta|\over |\xi_i-\eta_i|^q}  ~\lesssim~2^{-q|t_i-s_i|}.
 \eeq
For $t_i<s_i-\c$ and $q t_i\ge s_\jmath-\c$, the second equality in (\ref{norm 2 cases}) implies (\ref{iint norm est2}) again.
 
For $q t_i< s_\jmath-\c$, 
{\bf Remark 5.1} implies $|\xi_\jmath|\leq{1\over 2}|\eta_\jmath|$. Replace $i$ with $\jmath$ inside (\ref{int part Lambda_s}). We find 
 \bel{iint norm est1}   
{1+|\eta|\over |\xi_\jmath-\eta_\jmath|^q} ~\lesssim~2^{-(q-1)s_\jmath}~\lesssim~2^{-(q-1)s_i}.
 \eeq
 
 {\bf 2.}~Suppose $i\in\I_1\cap\J_2$. We have $|\xi_i|\sim 2^{t_i}$ and $|\eta_i|\lesssim 2^{s_i}$ where $|\eta|\sim|\eta_\jmath|\sim2^{s_\jmath}\sim2^{qs_i}$.

For $t_i>s_\ell-\c$, the first equality in (\ref{norm 2 cases}) remains valid which further implies (\ref{iint norm est2}).

For $t_i<s_\ell-\c$, {\bf Remark 5.1} implies 
$|\xi_\jmath|\leq{1\over 2}|\eta_\jmath|$.
Replace $i$ with $\jmath$ inside (\ref{int part Lambda_s}). We find (\ref{iint norm est1}). 
 \v
{\bf 3.}~Suppose $i\in\J_1\cap\I_2$. We have $|\xi_i|\lesssim 2^{t_i}$ where $|\xi|\sim|\xi_\imath|\sim2^{t_\imath}\sim2^{qt_i}$ and $|\eta_i|\sim 2^{s_i}$.

For $t_i>s_i+\c$, {\bf Remark 5.1} implies $|\eta_\imath|\leq{1\over 2}|\xi_\imath|$. 
Replace $i$ with $\imath$ inside (\ref{int part Lambda_s}). We find
 \bel{iint norm est3}   
{1+|\eta|\over |\xi_\imath-\eta_\imath|^q} ~\lesssim~2^{-(q-1)t_\imath}~\lesssim~2^{-(q-1)t_i}.
 \eeq
Consider $t_i<s_i-\c$ and $q t_i\ge s_\jmath-\c$. Note that the second equality in (\ref{norm 2 cases})
remains valid which implies (\ref{iint norm est2}). 

For  $q t_i< s_\jmath-\c$, {\bf Remark 5.1} implies $|\xi_\jmath|\leq{1\over 2}|\eta_\jmath|$. Replace $i$ with $\jmath$ inside (\ref{int part Lambda_s}). 
We find (\ref{iint norm est1}). 
\v
{\bf 4.}~Suppose $i\in\J_1\cap\J_2$. We have $|\xi|\sim 2^{q t_i}\sim2^{t_\imath}$ and $|\eta|\sim2^{q s_i}\sim2^{s_\jmath}$. 

For $t_i>s_i+\c$, we have $|\eta_\imath|\leq{1\over 2}|\xi_\imath|$. Replace $i$ with $\imath$ inside (\ref{int part Lambda_s}). We find (\ref{iint norm est3}).

For $t_i<s_i-\c$, we have $|\xi_\jmath|\leq{1\over 2}|\eta_\jmath|$.
Replace $i$ with $\jmath$ inside (\ref{int part Lambda_s}). We find (\ref{iint norm est1}). 
\v

In summary of step {\bf 1-4}, we integrate by parts $w.r.t$ $i$, $\imath$ and $\jmath$ depending on different cases. Moreover, the following estimates hold accordingly: 
\bel{term EST}
\begin{array}{cc}\ds
|\xi_i-\eta_i|~\sim~2^{t_i},~~~~ |\xi_i|~\sim~2^{t_i},~~~~ i\in\I_1,
\qquad
\hbox{or}\qquad |\xi_i-\eta_i|~\sim~2^{s_i},~~~~ |\eta_i|~\sim~2^{s_i},~~~~ i\in\I_2.
\\\\ \ds~~~
|\xi_\imath-\eta_\imath|~\sim~2^{t_\imath}, ~~~~ |\xi_\imath|~\sim~2^{t_\imath}
\qquad
\hbox{and}\qquad |\xi_\jmath-\eta_\jmath|~\sim~2^{s_\jmath},~~~~ |\eta_\jmath|~\sim~2^{s_\jmath}.
\end{array}
\eeq
Note that $|\xi_i|\gtrsim |\xi|^\rho$ because $|\xi|\sim|\xi_\imath|\sim2^{t_\imath}\lesssim 2^{q t_i}$ for $i\in\I_1$. Similarly, $|\eta_i|\gtrsim |\eta|^\rho$ because $|\eta|\sim|\eta_\jmath|\sim2^{t_\jmath}\lesssim 2^{q t_i}$ for $i\in\I_2$.

We carry out the integration by parts on every $i$ such that $t_i-s_i|>\c$. The resulting function, having an expression of $ \delta_\s(\eta) \prod\Delta_{y_i}\sigma(y,\eta)/ |\xi_i-\eta_i|^2\prod\Delta_{y_\imath}\sigma(y,\eta)/ |\xi_\imath-\eta_\imath|^2\prod\Delta_{y_\jmath}\sigma(y,\eta)/ |\xi_\jmath-\eta_\jmath|^2$
satisfies the differential inequality in (\ref{INEQ}). 

In particular, its norm is bounded by
\bel{norm t-s}
\C
\prod_{i=1}^n~2^{-\left({q-1\over q}\right)|t_i-s_i|}.
\eeq
We finish the proof by repeating the estimates from (\ref{Delta_k T dila}) to (\ref{Delta_t f Est Max}), with $\Omega_\t$ replaced by $\Omega_{\t\s}$ given in (\ref{Delta_t T Delta_s f}).   
\endproof

\section{Conclusion on the $\L^p$-boundedness}
\setcounter{equation}{0}
Let $\Delta_\t$ defined in (\ref{delta_t})-(\ref{Delta_t}). 
For $f\in\L^p(\R^\N)$ and $g\in\L^{p\over p-1}(\R^\N)$, we consider
\bel{Iden}
\int_{\R^\N} \T f(x) g(x) dx~=~\int_{\R^\N} \left\{\sum_\t \Delta_\t \T f(x)\right\} \left\{\sum_\s \Delta_\s g(x)\right\} dx.
\eeq
Note that  $\supp\deltaup_\t(\xi)\bar{\deltaup_\s(\xi)}\neq\emptyset$  only if $|\t-\s|\leq2$. 
Because of Plancherel theorem, our assertion for (\ref{Iden}) can be reduced to
\bel{Iden t}
\int_{\R^\N} \sum_\t \Delta_\t \T f(x) \Delta_\t g(x) dx.
\eeq
Let $h_i\in\Z, i=1,2,\ldots,n$. 
We define the $n$-tuple $(\t+\h)_i$ by simply replacing $t_i$ with $t_i+h_i$ respectively in (\ref{t_i dila}). 
For each $\t$ fixed, we write $f=\sum_{\h}\Delta_{\t+\h}f$. The expression in  (\ref{Iden t}) equals 
\bel{Iden est'}
 \sum_\h\int_{\R^\N} \sum_\t  \Delta_\t \T \Delta_{\t+\h}f(x) \Delta_\t g(x) dx.
\eeq
\begin{remark} We assume $\t, \t+\h\in\H$ satisfying (\ref{H}) in the summation of (\ref{Iden est'}).
\end{remark}
Note that $\Tilde{\delta}(\xi)=\sum_{\t\notin\H}\deltaup_\t(\xi)$ as shown in (\ref{Tilde delta rewrite}) has a compact support for  $|\xi|\lesssim 2^{\rhoup\over 1-\rhoup}$ for  $\xi\in\supp\Tilde{\deltaup}(\xi)$. 
The symbols of $\T\left(I-\sum_{\t+\h\notin\H}\Delta_{\t+\h}\right)$ and $\left[\left(I-\sum_{\t\notin\H}\Delta_\t\right)\T\right]^\ast$ respectively are $\Tilde{\delta}(\xi)\sigma(x,\xi)$ and $\bar{\Tilde{\delta}(\xi)\Lambda(x,\xi)}$ where $\Lambda(x,\xi)$ is given in (\ref{Lambda}). Both $\sigma\in\S_\rho, 0<\rho<1$ and $\Lambda(x,\xi)$ satisfy the differential inequality in (\ref{Ineq}).

From direct computation, we find
\bel{Ineq Tilde}
\begin{array}{lr}\ds
\left|\p_\xi^\alphaup\p_x^\betaup\Tilde{\delta}(\xi)\sigma(x,\xi)\right|~\leq~\C_{\alphaup~\betaup}~2^{\rho|\alpha|}~ \Bigg({1\over 1+|\xi|}\Bigg)^{|\alpha|}(1+|\xi|)^{\rhoup|\betaup|}
\end{array}
\eeq
for every multi-indices $\alphaup$, $\betaup$. The same is true for $\bar{\Tilde{\delta}(\xi)\Lambda(x,\xi)}$.

The regarding pseudo differential operators, as well as their adjoint operators,  are well known to be bounded on $\L^p$-spaces for $1<p<\infty$. 

By using Schwarz inequality and then H\"{o}lder inequality, we have 
\bel{projections}
\begin{array}{lr}\ds
\int_{\R^\N} \T f(x)g(x) dx
~\leq~\C  \sum_\h \int \left\{\sum_\t \Big(\Delta_{\t}\T \Delta_{\t+\h}f\Big)^2(x)\right\}^{1\over 2}\left\{\sum_\t\left(\Delta_\t g\right)^2(x)\right\}^{1\over 2} dx
\\\\ \ds ~~~~~~~~~~~~~~~~~~~~~~~~~~~~~
~\leq~\C \sum_\h \left\| \left\{\sum_{\t} \Big(\Delta_\t\T \Delta_{\t+\h} f\Big)^2\right\}^{1\over 2}\right\|_{\L^p(\R^\N)}\left\|\left\{\sum_\t\left(\Delta_\t g\right)^2\right\}^{1\over 2}\right\|_{\L^{p\over p-1}(\R^\N)}.
\end{array}
\eeq 
\begin{lemma}
Let $\Delta_\t$ defined in (\ref{Delta_t}). We have
\bel{Littlewood-Paley Ineq*}
\left\| \left\{\sum_\t\Big(\Delta_{\t}f\Big)^2\right\}^{1\over 2}\right\|_{\L^p(\R^\N)}
~\leq~\C_p~\left\|f\right\|_{\L^p(\R^\N)},\qquad 1<p<\infty.
\eeq
\end{lemma}
By taking the supremum of all $\|g\|_{\L^{p\over p-1}(\R^\N)}=1$ in (\ref{projections}) and applying {\bf Lemma 6.1}, we find 
\bel{sum}\left\|\T f\right\|_{\L^p(\R^\N)}~\leq~\C_p \sum_\h \left\| \left\{\sum_\t \Big(\Delta_\t\T\Delta_{\t+\h}f\Big)^2\right\}^{1\over 2}\right\|_{\L^p(\R^\N)},\qquad 1<p<\infty.
\eeq

Recall {\bf Lemma 5.1}. The norm of  $\Delta_\t \T f$ is bounded  by $\M f$ for every $\t\in\H$. 
By applying the vector-valued inequality of strong maximal function,  established by Fefferman and Stein \cite{C-F.S} and then using the Littlewood-Paley inequality in (\ref{Littlewood-Paley Ineq*}), we have
\bel{Max Littlewood-Paley Ineq}
\begin{array}{lr}\ds
  \left\| \left\{\sum_\t \Big(\Delta_\t\T\Delta_{\t+\h}f\Big)^2\right\}^{1\over 2}\right\|_{\L^p(\R^\N)}
~\leq~\C ~\left\| \left\{\sum_\t\Big(\M\Delta_{\t+\h}f\Big)^2\right\}^{1\over 2}\right\|_{\L^p(\R^\N)}
\\\\ \ds ~~~~~~~~~~~~~~~~~~~~~~~~~~~~~~~~~~~~~~~~~~~~~~~~
~\leq~\C_p~\left\| \left\{\sum_\t\Big(\Delta_{\t+\h}f\Big)^2\right\}^{1\over 2}\right\|_{\L^p(\R^\N)}
~\leq~\C_p~\left\|f\right\|_{\L^p(\R^\N)}
\end{array}
\eeq
for every $1<p<\infty$.

On the other hand, we have $\left(\Delta_\t \T \Delta_\s\right)^*f=\Delta^*_\s \T^* \Delta^*_\t f=0$. 
provided that $\supp\deltaup_\t(\xi)\bar{\deltaup_\s(\xi)}=\emptyset$ if $|\t-\s|>2$.  
Given $\h$ fixed, by applying Cotlar-Stein lemma  together with {\bf Lemma 6.1}, we have 
\bel{L2 bound}
\left\|\left\{\sum_\t\Big(\Delta_{\t} \T \Delta_{\t+\h} f\Big)^2\right\}^{1\over 2}\right\|_{\L^2(\R^\N)}~\leq~\C_\rho~\prod_{i=1}^n 2^{-(1-\rhoup) |h_i|}\left\|f\right\|_{\L^2(\R^\N)}.
\eeq
From  (\ref{Max Littlewood-Paley Ineq}) and (\ref{L2 bound}),    Riesz-Thorin interpolation theorem implies
\bel{Lpbound}
 \left\| \left\{\sum_\t \Big(\Delta_\t\T\Delta_{\t+\h}f\Big)^2\right\}^{1\over 2}\right\|_{\L^p(\R^\N)}~\leq~\C_{\rho~p}~\prod_{i=1}^n 2^{-\ve |h_i|}\left\|f\right\|_{\L^p(\R^\N)}
\eeq
for $p\ne2, 1<p<\infty$ where $\ve=\ve(\rhoup,p)>0$. By summing over all the $h_i, i=1,2,\ldots,n$   in (\ref{sum}),  we obtain the desired result.

\section{Proof of the Littlewood-Paley inequality}
\setcounter{equation}{0}
In order to prove {\bf Lemma 6.1}, we momentarily consider $n=2$. 

Let $(x,y)\in\R^{\N_1}\times\R^{\N_2}$ and 
 $(\xi,\eta)\in\R^{\N_1}\times\R^{\N_2}$.  
 
 Recall $\phi$ defined in (\ref{delta_t}).  For $0<\rhoup<1$ and $\delta>0$, we write
\bel{Phi_delta'}
\Hat{\Phi}_\delta(\xi,\eta)~=~\phi\left(\delta \xi, \delta^{1\over\rhoup} \eta\right).
\eeq
Note that
\bel{cancellation Phi}
\iint_{\R^{\N_1}\times\R^{\N_2}} \Phi_\delta(x,y)dxdy~=~\phi(0,0)~=~0,\qquad \delta>0.
\eeq
Moreover,  $\Phi(x,y)\doteq\Phi_1(x,y)$ is smooth, bounded and decays rapidly as $|(x,y)|\mt\infty$. 

The regarding square function is defined by 
\bel{Square function}
\Sz_\Phi f(x,y)~=~\left\{\int_{0}^\infty \Big( f\ast\Phi_{\delta}\Big)^2(x,y) {d\delta\over \delta}\right\}^{1\over 2}.
\eeq
We aim to show
\bel{Square function RESULT}
\left\| \Sz_\Phi f\right\|_{\L^p(\R^\N)}~\leq~\C_p~\left\| f\right\|_{\L^p(\R^\N)},\qquad 1<p<\infty.
\eeq
The kernel of $\Sz_\Phi$, denoted by 
\bel{Phi_delta}
\Omega(x,y)~=~ \Phi_\delta(x,y)~=~\left({1\over \delta}\right)^{\N_1+{\N_2/\rhoup}}\Phi\left({x\over \delta},{y\over \delta^{1/\rhoup}}\right)
 \eeq
has a Hilbert space valued norm 
\bel{Omega^delta norm}
\left|\Omega(x,y)\right|_{H}~=~\left\{\int_{0}^\infty \left|\Phi_{\delta}(x,y)\right|^2 {d\delta\over \delta}\right\}^{1\over 2}.
\eeq

\begin{lemma}
\bel{Omega^delta Est}
\left|\p_x^\alphaup\p_y^\betaup\Omega(x,y)\right|_{H}~\leq~\C_{\alphaup~\betaup}~\left({1\over |x|+|y|^\rhoup}\right)^{\N_1+|\alphaup|}\Bigg({1\over |y|+|x|^{1/ \rhoup}}\Bigg)^{\N_2+|\betaup|}
\eeq
for every multi-indices $\alphaup$, $\betaup$.
\end{lemma}

{\bf Proof}~ 
{\bf 1.} Suppose $|x|\leq\delta$ and $|y|\leq\delta^{1/\rhoup}$. Because $\Phi$ is smooth and  bounded, (\ref{Phi_delta}) implies 
\bel{Omega^delta est1 x}
\begin{array}{lr}\ds
 \left|\Omega(x,y)\right|_{H}~\leq~\C\left\{\int_{|x|}^\infty \left({1\over |x|}\right)^{2\N_1+2\N_2/\rhoup-1} {d\delta\over \delta^2}\right\}^{1\over 2}
 ~\leq~\C \left({1\over |x|}\right)^{\N_1+\N_2/\rhoup}
 \end{array}
 \eeq
and 
\bel{Omega^delta est1 y}
\begin{array}{lr}\ds
 \left|\Omega(x,y)\right|_{H}~\leq~\C\left\{\int_{|y|^\rhoup}^\infty \left({1\over |y|}\right)^{2\rhoup\N_1+2\N_2-1} {d\delta\over \delta^{1+1/\rhoup}}\right\}^{1\over 2}
 ~\leq~\C \left({1\over |y|}\right)^{\rhoup\N_1+\N_2}.
 \end{array}
 \eeq
 
{\bf 2.} Suppose $|x|>\delta$ and $|y|\leq \delta^{1/\rhoup}$. We necessarily have $|x|>|y|^\rhoup$. Because $\Phi(x,y)$ decays rapidly as $|(x,y)|\mt\infty$, we have 
\bel{Omega^delta est2 x}
\begin{array}{lr}\ds
 \left|\Omega(x,y)\right|_{H}~\leq~\C\left\{\int_0^{|x|} \left({1\over |x|}\right)^{2\N_1+2\N_2/\rhoup+2}\delta d\delta\right\}^{1\over 2}
 ~\leq~\C \left({1\over |x|}\right)^{\N_1+\N_2/\rhoup}.
 \end{array}
 \eeq
 On the other hand, the estimate in (\ref{Omega^delta est1 y}) shows 
 \bel{Omega^delta est2 y}
 \left|\Omega(x,y)\right|_{H}~\leq~\C \left({1\over |y|}\right)^{\rhoup\N_1+\N_2}.
 \eeq
 
{\bf 3.} Suppose $|x|\leq\delta$ and $|y|>\delta^{1/\rhoup}$. We necessarily have $|x|<|y|^\rhoup$. The estimate in (\ref{Omega^delta est1 x}) shows 
\bel{Omega^delta est3 x}
 \left|\Omega(x,y)\right|_{H}~\leq~\C \left({1\over |x|}\right)^{\N_1+\N_2/\rhoup}.
 \eeq
On the other hand, because $\Phi(x,y)$ decays rapidly as $|(x,y)|\mt\infty$, we have 
\bel{Omega^delta est3 y}
\begin{array}{lr}\ds
 \left|\Omega(x,y)\right|_{H}~\leq~\C\left\{\int_0^{|y|^\rhoup} \left({1\over |y|}\right)^{2\rhoup\N_1+2\N_2+2}\delta^{2/\rhoup-1} d\delta\right\}^{1\over 2}
 ~\leq~\C \left({1\over |y|}\right)^{\rhoup\N_1+\N_2}.
 \end{array}
 \eeq
 
{\bf 4.} Suppose $|x|>\delta$ and $|y|>\delta^{1/\rhoup}$. From (\ref{Omega^delta est2 x}) and (\ref{Omega^delta est3 y}), we find
\bel{Omega^delta est4 }
 \left|\Omega(x,y)\right|_{H}~\leq~\C\left({1\over |x|}\right)^{\N_1+\N_2/\rhoup},\qquad  \left|\Omega(x,y)\right|_{H}~\leq~\C \left({1\over |y|}\right)^{\rhoup\N_1+\N_2}.
 \eeq
All together, we obtain 
\bel{Omega^delta est}
\left|\Omega(x,y)\right|_{H}~\leq~\C\left({1\over |x|+|y|^\rhoup}\right)^{\N_1}\Bigg({1\over |y|+|x|^{1/ \rhoup}}\Bigg)^{\N_2}.
\eeq
Observe that every $\p_x$ acting on $\Phi_\delta$ gains a constant multiple of $\delta^{-1}$ and every $\p_y$ acting on $\Phi_\delta$ gains a constant multiple of $\delta^{-1/\rhoup}$ respectively. By carrying out the same estimates as above, we prove the differential inequality in (\ref{Omega^delta Est}).
\endproof

Define 
\bel{dist fun}
\lambda(x,y)~=~\max\Big\{|x|,|y|^\rhoup\Big\}
\eeq
and the non-isotropic ball
\bel{non-iso ball}
\B_{\lambda,\delta}~=~\Big\{(x,y)\in\R^{\N_1}\times\R^{\N_2}~\colon~\lambda(x,y)\leq\delta\Big\}.
\eeq
Let $(u,v)\in\B_{\lambda,\delta}$ and $(x,y)\in {^c\B}_{\lambda, r\delta}$ for some $r>1$. By adjusting its value, we can have
\bel{rho u v x y}
\lambda(x,y)~>~r\delta,\qquad \lambda(x-u,y-v)~>~\delta.
\eeq
The differential inequality in (\ref{Omega^delta Est}) implies 
\bel{Omega^delta difference}
\begin{array}{lr}\ds
\left|\Omega(x-u,y-v)-\Omega(x,y)\right|_{H}
\\\\ \ds 
~\leq~\C~{|u|\over\lambda(x,y)}\left({1\over\lambda(x,y)}\right)^{\N_1+\N_2/\rhoup}+\C~{|v|\over\lambda(x,y)^{1/ \rhoup}}\left({1\over\lambda(x,y)}\right)^{\N_1+\N_2/\rhoup}
\\\\ \ds 
~\leq~\C~{\lambda(u,v)\over\lambda(x,y)}\left({1\over\lambda(x,y)}\right)^{\N_1+\N_2/\rhoup}.
\end{array}
\eeq
From (\ref{rho u v x y}) and (\ref{Omega^delta difference}), we have
 \bel{Omega differen Est}
\begin{array}{lr}\ds
\sum_{k=0}^\infty \iint_{ \B_{\lambda,2^{k+1}r\delta}\setminus\B_{\lambda,2^kr\delta}} \left|\Omega(x-u,y-v)-\Omega(x,y)\right|_{H} dxdy
\\\\ \ds
~\leq~\C\sum_{k=0}^\infty \iint_{ \B_{\lambda,2^{k+1}r\delta}\setminus\B_{\lambda,2^kr\delta}} {\lambda(u,v)\over\lambda(x,y)}\left({1\over\lambda(x,y)}\right)^{\N_1+\N_2/\rhoup}dxdy
\\\\ \ds
~\leq~\C~2^{\N_1+\N_2/\rhoup} \sum_{k=0}^\infty  2^{-k},\qquad (u,v)\in\B_{\lambda,\delta}.
\end{array}
\eeq
From (\ref{Omega differen Est}), for every $(u,v)\in\B_{\lambda,\delta}$, we conclude
\bel{Omega^delta differen Est}
\iint_{{^c}\B_{\lambda,r\delta}}\left|\Omega(x-u,y-v)-\Omega(x,y)\right|_{H} dxdy~\leq~\C.
\eeq
This is a well known condition for $f\ast\Omega$ satisfying the  weak type $(1,1)$-estimate with the Hilbert space valued norm in  (\ref{Phi_delta})-(\ref{Omega^delta norm}).  More regarding details can be found in chapter I of  Stein \cite{S}. 

On the other hand, by applying Plancherel theorem, we have
\bel{S_phi L^2}
\begin{array}{lr}\ds
\left\|\Sz_\Phi f\right\|_{\L^2(\R^{\N_1}\times\R^{\N_2})}^2~=~\iint_{\R^{\N_1}\times\R^{\N_2}}\left\{\int_0^\infty |\Hat{f}(\xi,\eta)|^2\left|\Hat{\Phi}\left(\delta\xi,\delta^{1\over \rhoup}\eta\right)\right|^2{d\delta\over \delta}\right\}d\xi d \eta
\\\\ \ds ~~~~~~~~~~~~~~~~~~~~~~~~~~~~~~
~\leq~\left\{\sup_{(\xi,\eta)}\int_0^\infty \phi^2\left(\delta\xi,\delta^{1\over \rhoup}\eta\right){d\delta\over \delta}\right\}\iint_{\R^{\N_1}\times\R^{\N_2}} |\Hat{f}(\xi,\eta)|^2d\xi d\eta
\\\\ \ds ~~~~~~~~~~~~~~~~~~~~~~~~~~~~~~
~\leq~\C~\left\|f\right\|_{\L^2(\R^{\N_1}\times\R^{\N_2})}^2.
\end{array}
\eeq
From the weak type $(1,1)$-estimate and (\ref{S_phi L^2}), we conclude the $\L^p$-result in (\ref{Square function RESULT}) by using Marcinkewicz interpolation theorem and a standard argument of duality.

Consider 
$\Hat{\Psi}_\delta(\xi,\eta)=\phi\left(\delta^{1\over\rhoup} \xi, \delta \eta\right)$. Our estimates proving {\bf Lemma 7.1} remain  valid for $\Omega\doteq\Psi_\delta$, with $x$ and $y$ switched in roles. Denote $F=f\ast\Phi_\delta$. By applying  the $\L^p$-regularity theorem, stated as Theorem 3 in chapter I of  Stein \cite{S} for space valued functions,  we have 
\bel{Iteration est}
\begin{array}{lr}\ds
\iint_{\R^{\N_1}\times\R^{\N_2}}\Big(\Sz_{\Phi\ast\Psi} f\Big)^p(x,y)dxdy~=~\iint_{\R^{\N_1}\times\R^{\N_2}} \left\{\int_{0}^\infty \left| F \ast\Psi_\delta(x,y)\right|_H^2 {d\delta\over\delta}\right\}^{p\over 2}dxdy
\\\\ ~~~~\qquad\qquad\qquad\qquad\qquad~~~~~~~~\ds
~\leq~\C_p~\iint_{\R^{\N_1}\times\R^{\N_2}} \left| F(x,y)\right|_H^p dxdy
\\\\ ~~~~\qquad\qquad\qquad\qquad\qquad~~~~~~~~\ds 
~=~\C_p~\iint_{\R^{\N_1}\times\R^{\N_2}} \left\{\int_0^\infty \left|f\ast\Phi_\delta(x,y)\right|^2{d\delta\over \delta}\right\}^{p\over 2} dxdy
\\\\ ~~~~\qquad\qquad\qquad\qquad\qquad~~~ ~~~~~\ds 
~\leq~\C_p~\iint_{\R^{\N_1}\times\R^{\N_2}} \left|f(x,y)\right|^p dxdy,\qquad 1<p<\infty.
\end{array}
\eeq
Lastly, it is clear that the above iteration argument works for any number of parameter. Replace $\N_1, \N_2$ with $\N_i, \N-\N_i$ for every $i=1,2,\ldots,n$ in (\ref{Phi_delta'})-(\ref{S_phi L^2}). By carrying out an $n$-parameter analogue of (\ref{Iteration est}), we conclude the $\L^p$-boundedness of the corresponding square function. This is equivalent to our desired Littlewood-Paley inequality in (\ref{Littlewood-Paley Ineq*}).

\v

{\bf Acknowledgement

~~~~~~~ I am deeply grateful to my advisor Elias M. Stein for those stimulating talks and unforgettable lectures.}

\end{document}